# New Confocal Hyperbola-based Ellipse Fitting with Applications to Estimating Parameters of Mechanical Pipes from Point Clouds

Reza Maalek[1,*], Derek D. Lichti[2]

[1] Endowed Professor of Digital Engineering and Construction, Institute for Technology and Management in Construction (TMB), Karlsruhe Institute of Technology (KIT). Address: Am Fasanengarten, Karlsruhe, Germany, 76131. Email: reza.maalek@kit.edu (corresponding author)

[2] Professor of Geomatics Engineering, Department of Geomatics Engineering of the Schulich School of Engineering, University of Calgary. Address: 2500 University Dr. NW, Calgary, Canada, T2N 1N4. Email: ddlichti@ucalgary.ca

* Indicates corresponding author

**Abstract:** This manuscript presents a new method for fitting ellipses to two-dimensional data using the confocal hyperbola approximation to the geometric distance of points to ellipses. The proposed method was evaluated and compared to established methods on simulated and real-world datasets. First, it was revealed that the confocal hyperbola distance considerably outperforms other distance approximations such as algebraic and Sampson. Next, the proposed ellipse fitting method was compared with five reliable and established methods proposed by Halir, Taubin, Kanatani, Ahn and Szpak. The performance of each method as a function of rotation, aspect ratio, noise, and arc-length were examined. It was observed that the proposed ellipse fitting method achieved almost identical results (and in some cases better) than the gold standard geometric method of Ahn and outperformed the remaining methods in all simulation experiments. Finally, the proposed method outperformed the considered ellipse fitting methods in estimating the geometric parameters of cylindrical mechanical pipes from point clouds. The results of the experiments show that the confocal hyperbola is an excellent approximation to the true geometric distance and produces reliable and accurate ellipse fitting in practical settings.

## 1. Introduction: applications of ellipse fitting

Points following elliptic patterns commonly and naturally occur in many real-world datasets. In digital images, the edge points of spheres and circles form ellipses due to projective transformation. In three-dimensional point clouds, the intersection of spheres, cones, and cylinders with a plane also delivers points following an elliptic pattern. Reliable and accurate ellipse fitting to data is, hence, commonly an important step towards solving various real-world applications pertaining to shapes such as circles, cylinders and spheres, including, galaxy classification [1], traffic sign detection [2], phase estimation in interferograms [3], cell segmentation in fluorescence microscopy images [4], demodulation of interferometer signals [5], terrestrial laser scanner (TLS) registration [6], optical instrument calibration [7], geometric as-built building information modeling (BIM) of pipes [8], cylinder fitting of industrial shapes [9], and cylindrical point cloud classification [10].

To fit an ellipse to a set of points, the ellipse parameters minimizing the sum of the squared distances between the points and an ellipse are estimated (least squares method). Ideally, it is desirable to minimize the function representing the orthogonal distances of points to an ellipse (commonly referred to as the geometric distance). Finding the orthogonal contact point to an ellipse is, however, not an easy task, and most available analytical methods are either numerically unstable, or iterative without guarantee of convergence [11], and consequentially time consuming [12]. The numerically stable methods, such as that proposed by Chernov [11] are not analytically closed-form, and hence, cannot be used as a minimization cost function. As a result, researchers have sought to minimize alternative and less complicated distance functions, such as algebraic [13–16], and Sampson [12,17–20]. The methods minimizing algebraic and Sampson distances, however, share some undesirable properties. Algebraic methods are biased at locations of high curvature (see Section 2). The Sampson distance is not continuous in $\mathcal{R}^2$ (see Sections 2 and 5.1.1 Figure 4) and the methods are only suited for data with "moderate" noise levels [12]. The methods minimizing the algebraic and Sampson distances also require either additional constraints [15] or barrier terms [12] during the minimization of the objective function to guarantee an ellipse. Furthermore, Sampson and algebraic distances do not effectively predict the geometric distance and as will be revealed in our experiments (see Section 5.1.1 and Figure 4), diverge from the true geometric distance as the point moves farther from the ellipse. This is most likely the reason they are only suited for "moderate" noise levels.

To overcome the discussed limitations associated with current ellipse fitting methods, this manuscript presents a new ellipse fitting method that minimizes the sum of the squared confocal hyperbola distances of points to an ellipse. Even though the confocal hyperbola distance approximation was introduced in Rosin [21] with promising geometric properties, its effectiveness as an ellipse fitting method has not yet been examined in the literature. The goals of this manuscript are to: (i) provide a closed-form solution to the confocal hyperbola distance function; (ii) introduce a new algorithm to minimize the sum of the squared confocal hyperbola distances; and (iii) thoroughly investigate its effectiveness in fitting ellipses to two-dimensional (2D) data, acquired from simulated image edge points as well as point clouds of cylindrical mechanical pipes.

The remainder of the manuscript is structured as follows. In Section 2, an introduction to current ellipse fitting methods, along with their advantages and limitations is provided. Section 3 includes the detailed explanation of the proposed ellipse fitting algorithm. Section 4 outlines the metrics for validation of the proposed method along with the design and configurations of the various simulation-based and real-world experiments, carried out to assess the applicability of the proposed method compared to state-of-the-art ellipse fitting methods. Sections 5 and 6 detail the results obtained by the simulated, and real-world data in the designed configurations presented in Section 4, respectively. Section 7 discusses the summary of the major findings of this study.

## 2. Background: introduction to ellipse fitting

The equation of an ellipse in general form is represented as:

$$(\frac{(x-x_c)\cos\theta+(y-y_c)\sin\theta}{a})^2+(\frac{-(x-x_c)\sin\theta+(y-y_c)\cos\theta}{b})^2=1 \quad (1)$$

where $(x, y)$ are the point coordinates, $\rho = (x_c, y_c, a, b, \theta)$ is the geometric parameter vector of an ellipse, which includes the coordinates of the center, semi major length, semi minor length and rotation angle of the ellipse, respectively. Equation (1) can be reduced to general algebraic form as proposed by Paton [13], which is commonly referred to as the algebraic distance:

$$d_{alg}(x,y)=Ax^2+Bxy+Cy^2+Dx+Ey+F=\tau\, P^T=0 \quad (2)$$

where $d_{alg}(x,y)$ is the algebraic distance of point $(x,y)$ to the conic, $\tau=(A,B,C,D,E,F)$ is the vector of algebraic parameters, $P=(x^2,xy,y^2,x,y,1)$ is the 6-dimensional coordinate set (commonly referred to as the design vector, or for more than one point, design matrix), and $(.)^T$ denotes the transpose function. To find a best fit conic to $N$ data points, the parameters providing the minimum sum of the squared algebraic distances are required, i.e.:

$$\min_{\tau}\sum_{i=1}^{N}\left(d_{alg}(x_i,y_i)\right)^2=\min_{\tau}\tau(P^TP)\tau^T=\min_{\tau}\tau S\tau^T \quad (3)$$

where $S=P^TP$ is referred to as the scatter matrix. The problem imposed by equation 3 is a linear least squares problem subject to some parameter constraint to prevent the trivial solution, $\tau = 0$. For instance, Paton [13] proposed to minimize equation 3, subject to $||\tau|| = 1$, where $||.||$ denotes the L2-norm. This problem is now reduced to fitting a hyper plane to a six-dimensional point coordinate set $P_i=({x_i}^2,x_iy_i,{y_i}^2,x_i,y_i,1)$, which is the eigenvector corresponding to the smallest eigenvalue of the 6×6 scatter matrix, $S$. Since this algebraic minimization ignores the correlations between the algebraic parameters, it cannot guarantee an ellipse. In fact, the algebraic parameters must satisfy both of the following conditions to produce an ellipse:

$$Ellipse = \begin{cases} D = \det\left(\begin{bmatrix} A & \frac{B}{2} & \frac{D}{2} \\ \frac{B}{2} & C & \frac{E}{2} \\ \frac{D}{2} & \frac{E}{2} & F \end{bmatrix}\right) & , D \neq 0 \\ \Delta = B^2 - 4AC & , \Delta < 0 \end{cases} \quad (4)$$

where $\det(.)$ denotes the determinant function. To guarantee an ellipse fit, Fitzgibbon [15] proposed to solve equation 3 subject to $\Delta = -1$. By introducing a constraint matrix ($M$) for this condition, the problem was reduced to a generalized eigenvalue problem ($S\tau^T = \lambda M\tau^T$). Halir [16] observed that the constraint matrix of Fitzgibbon can provide numerically unstable results, especially when the points are exactly on the ellipse. Hence, they proposed a numerically stable method to solve Fitzgibbon's constraint.

To provide some perspective on the meaning of the algebraic distance, Bookstein [14] showed that the algebraic distance is proportional to:

$$d_{alg}(x_i, y_i) \propto \frac{d_i^2 - c_i^2}{c_i^2} = \frac{l_i}{c_i}\left(\frac{l_i}{c_i} + 2\right) \quad (5)$$

where $d_i$ is the distance from the point to the center, and $c_i$ and $l_i$ are shown in Figure 1a. The formulation of Bookstein shows that for the same $l_i$ across all data points (e.g. same noise level), the algebraic distance of the point to the curve is smaller when $c_i$ is larger (i.e. around the major axis). In other words, at locations of high curvature, the algebraic distance under-estimates the distance of the point to the curve. This is referred to as a high curvature bias problem imposed by equation 3 [22]. To solve the curvature bias problem, Agin [17] proposed the minimization of the following distance function:

$$d_{Sampson}(x_i, y_i) = \frac{d_{alg}(x_i, y_i)}{\|\nabla d_{alg}(x_i, y_i)\|} \quad (6)$$

where $\nabla$ denotes the gradient function, and $d_{Sampson}(x, y)$ is commonly referred to as the Sampson distance of points $(x, y)$ to the ellipse (Sampson 1982). Minimization of the sum of the squared Sampson distances is, however, a non-linear problem. As a result, Agin [17] proposed to approximate equation 6 subject to $\sum_{i=1}^{N} \|\nabla d_{alg}(x_i, y_i)\|^2 = 1$, which reduces to a linear generalized eigenvalue problem. Later, Taubin [19] independently derived this exact formulation, which is now commonly referred to as the Taubin's method in the literature. Instead of solving the actual non-linear minimization problem of equation 6, Sampson [18] proposed to adopt a reweighting strategy (referred to as the gradient weighted method) using the distance of equation 6 to iteratively improve Bookstein's algebraic conic fitting method. The approximate solution to equation 6 proposed by Taubin and Agin provides a curvature bias correction of the algebraic distance up to the first order. Kanatani [23] went one step further and corrected the algebraic distance approximation up to the second order of the algebraic distance. Their method also reduces to a generalized eigenvalue minimization problem similar to that of Taubin and Agin. The approximate solutions to equation 6 using Taubin's and Kanatani's methods, however, do not necessarily satisfy the conditions of equation 4, and hence, cannot guarantee an ellipse. To this end, Szpak [12] proposed to incorporate a barrier term in the minimization of the Sampson distance, which helps limit the solution space to guarantee an ellipse. Their formulation, however, cannot be reduced to a linear least squares problem, and hence it was solved using the Levenberg–Marquardt algorithm.

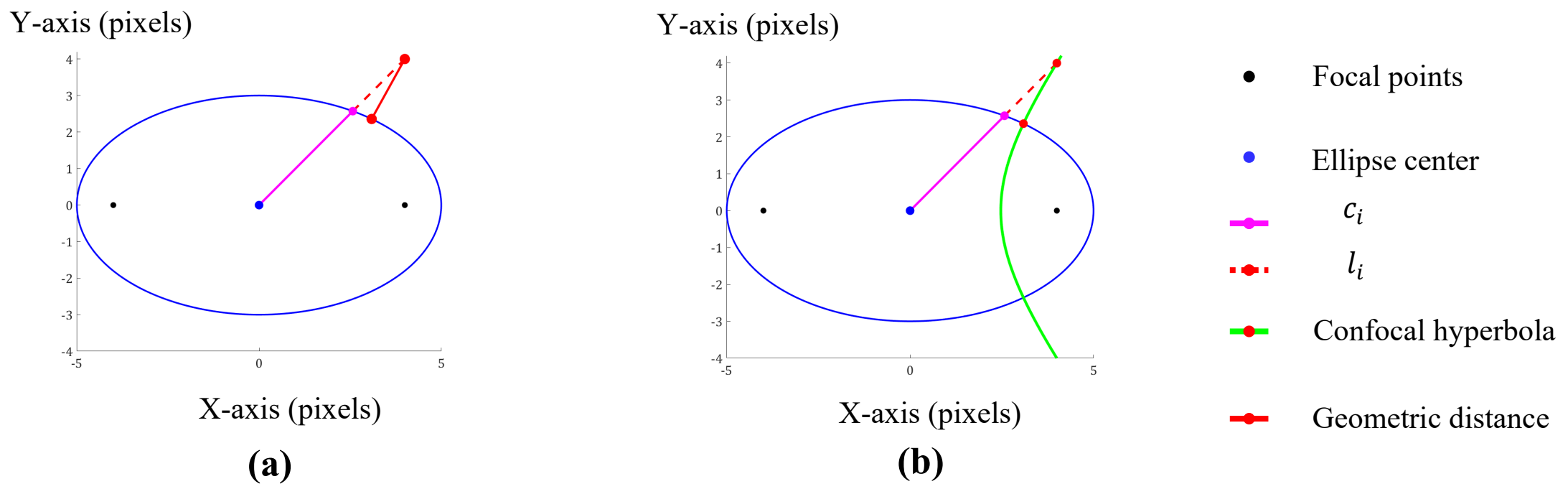

**Figure 1:** Schematic representation of distance function for one point to an ellipse with geometric parameter vector $\rho = (0,0,5,3,0)$: a) algebraic and geometric distances; b) confocal hyperbola

The Sampson distance, although providing a good approximation to reduce the impact of curvature biases imposed by the algebraic distance function, is undetermined when the point coordinate is on the center of the ellipse (i.e. equation 6 goes to infinity as the point approaches the center of the ellipse). In addition, it bears no direct relationship to the true geometric distance, especially as the point moves farther from the ellipse (see visual analysis of Section 5.1.1). Harker and O'Leary [24] pointed out that the Sampson distance is merely the geometric distance of a point to the first order approximation of the algebraic distance function and should not be mistaken with the first order approximation of the geometric distance from the point to the conic. They provided a closed-form solution to the first order approximation of the geometric distance, which was, in fact, also found to be a weighted algebraic distance. However, as pointed out by Uteshev and Goncharova [25], Harker and O'Leary's solution can display singularities in elongated ellipses for points around the major axis. In addition, as will be revealed later (see Figure 4), Haker and O'Leary's distance is undetermined (approaches infinity) for points at the focal points of an ellipse.

The ideal geometric minimization cost function is to minimize the sum of the squared distances represented by the orthogonal projections of the points onto the ellipse (i.e. the geometric distance; Figure 1a; red solid line).

Achieving a numerically stable, closed-form analytical solution to the geometric distance is, however, not an easy task [11]. For instance, Safaee-Rad [26] proposed an exact solution to the geometric distance that reduces to finding the roots of a quartic equation. However, analytical methods for finding the roots of a quartic equation are known to be numerically unstable. To this end, Ahn [27] proposed to divide the solution provided by Safaee-Rad into two functions that were iteratively solved using the generalized Gauss-Newton method. However, as pointed by Chernov [11], the method does not guarantee convergence and numerically failed about 3% of the time. Currently, the most accurate and numerically stable process for finding the orthogonal projection of points onto conics is proposed by Chernov [11]. The method converts the orthogonal projection problem to finding one real generalized eigenvalue between the ellipse and its auxiliary conic (which turns out to be a hyperbola; see Figure 1 of [11]). Using one real eigenvalue root, a degenerative conic, comprised of two lines, is produced. The closest orthogonal projection of a point to the ellipse is the point of intersection between the two lines and the ellipse (up to four points), whose distance to the point of interest is the least. However, this elegant method is still not analytically closed-form and, hence, cannot be used as a minimization cost function to estimate the parameters of the best fit ellipse.

Another approximation to the geometric distance is the confocal hyperbola method used in Rosin [21], see Figure 1b. This distance approximation relies on the fact that confocal ellipses and hyperbolas are mutually orthogonal at the intersecting point. Since the confocal hyperbola passing through a given point is nearly linear at the point of intersection, it is a good approximation to the geometric distance of a point to the ellipse. It is worth mentioning that the auxiliary conic to an ellipse at a specific point, as shown by Chernov [11], is also a hyperbola. Hence, the confocal hyperbola distance is speculated to possess similar geometric properties to the true geometric distance of a point to the ellipse. In fact, Rosin had compared 16 different geometric distance approximation methods, including algebraic, Sampson, Harker and O'Leary, and the confocal hyperbola, in terms of properties such as linearity, curvature bias, asymmetry and general goodness [21,28,29]. Overall, the confocal hyperbola was the best in every category, especially for larger noise levels (see Table 3 of [21,28,29]). In Section 5.1, we will demonstrate that the confocal hyperbola method provides almost identical results to that of Chernov [11], which is the current state-of the-art. Although Rosin provided a process to derive the confocal hyperbola distance, it was not presented in general closed-form, imposed singularities for circles, and required the examination of four solutions. As a result, up until now, the confocal hyperbola distance had not been used to fit ellipses to data. In Section 3, we will expand on the formulation provided by Rosin and present an analytically closed-form solution to this geometric distance approximation, which can then be minimized to obtain the best fit ellipse parameters.

## 3. Ellipse fitting using confocal hyperbola distance

Various ellipse fitting methods along with their limitations were introduced in Section 2. Here, we will formulate a new ellipse fitting method based on the minimization of the confocal hyperbola distance. To produce a minimization objective function, we shall first derive a closed-form solution to the confocal hyperbola distance by determining the point of intersection between the ellipse and the confocal hyperbola passing through a given point $(x, y)$. Using a translation and rotation, equation 1 can be rewritten as follows:

$$\frac{X^2}{{a_e}^2}+\frac{Y^2}{{b_e}^2}=1 \text{ with } X=(x-x_c)\cos\theta+(y-y_c)\sin\theta \text{ and } Y=-(x-x_c)\sin\theta+(y-y_c)\cos\theta \quad (7)$$

where $a_e$ and $b_e$ are the semi-major and semi-minor lengths of the ellipse, respectively. In this new coordinate system, the problem is now reduced to finding the point of intersection between the ellipse of equation (7) and the confocal hyperbola passing through point $(X, Y)$ with equation:

$$\frac{X^2}{{a_h}^2}-\frac{Y^2}{{b_h}^2}=1 \quad (8)$$

where $a_h$ and $b_h$ are the semi-major and semi-minor lengths of the hyperbola, respectively. Since the ellipse and hyperbola are confocal, the following relationships are valid:

$$f^2=a_e^2-b_e^2=a_h^2+b_h^2 \quad (9)$$

where $f$ is the distance of one of the foci to the center. By substituting equation 9 into 8, the following quadratic function is derived:

$$A^2-(X^2+Y^2+f^2)A+X^2f^2=0 \text{ with } A={a_h}^2 \quad (10)$$

Solving equation 10 for $A$ provides two solutions; however, only the following can ensure both $A$ and $f^2-A$ are positive (this is of course a necessary condition to produce real values for $a_h$ and $b_h$):

$$\begin{cases} A={a_h}^2=\frac{X^2+Y^2+f^2-\sqrt{(X^2+Y^2+f^2)^2-4X^2f^2}}{2} \\ f^2-A={b_h}^2=\frac{\sqrt{(X^2+Y^2-f^2)^2+4Y^2f^2}-(X^2+Y^2-f^2)}{2} \end{cases} \quad (11)$$

Substituting equation (11) into (8) and solving simultaneously with (7) will provide the following simplified solutions to the orthogonal intersecting point $(X_I, Y_I)$:

$$|X_I|=\frac{a_e a_h}{f}\ ,\ |Y_I|=\frac{b_e b_h}{f} \quad (12)$$

Since four solutions exists $(\pm X_I, \pm Y_I)$, Rosin [21] proposed to calculate the distances between the point and the four solutions of equation 12 and select the solution with the minimum distance. However, it can be observed that the minimum distance must belong to the solution that is in the same quadrant as the point. Since $a_e, a_h, b_e, b_h$ and $f$ are all positive, the correct solution can be formulated by simply restricting the point coordinate $(X, Y)$ to the first quadrant using an absolute value function. The confocal hyperbola distance can now be defined in closed form as follows:

$$D_h(X,Y) = \left\| \begin{bmatrix} D_{hX}(X,Y) \\ D_{hY}(X,Y) \end{bmatrix} \right\| = \left\| \begin{bmatrix} |X| - |X_I| \\ |Y| - |Y_I| \end{bmatrix} \right\| \tag{13}$$

where $|.|$ is the absolute value function, $\begin{bmatrix} D_{hX}(X,Y) \\ D_{hY}(X,Y) \end{bmatrix}$ is a 2D vector of distances in the translated and rotated $X - Y$ coordinate plane (see equation 7), and $\|.\|$ denotes the L2-norm. In its current form, equation 12 is numerically singular when $a_e = b_e$ (i.e. $f = 0$, when the data follow a circular pattern). However, if the formulations of $a_h$ and $b_h$ in equation 11 are rearranged, we will have:

$$\begin{cases} {a_h}^2 = \frac{X^2+Y^2+f^2-\sqrt{(X^2+Y^2+f^2)^2-4X^2f^2}}{2} \cdot \frac{X^2+Y^2+f^2+\sqrt{(X^2+Y^2+f^2)^2-4X^2f^2}}{X^2+Y^2+f^2+\sqrt{(X^2+Y^2+f^2)^2-4X^2f^2}} = \frac{2X^2f^2}{X^2+Y^2+f^2+\sqrt{(X^2+Y^2+f^2)^2-4X^2f^2}} \\ {b_h}^2 = f^2 - {a_h}^2 = f^2 - \frac{2X^2f^2}{X^2+Y^2+f^2+\sqrt{(X^2+Y^2+f^2)^2-4X^2f^2}} \end{cases} \tag{14}$$

therefore, the closed formula for the orthogonal intersecting point of equation 12 will now become:

$$\begin{cases} |X_I| = \frac{a_e a_h}{f} = \frac{a_e \cdot |X|}{\sqrt{\frac{X^2+Y^2+f^2+\sqrt{(X^2+Y^2+f^2)^2-4X^2f^2}}{2}}} \\ |Y_I| = \frac{b_e b_h}{f} = \frac{b_e}{a_e}\sqrt{a_e^2 - X_I^2} \end{cases} \xrightarrow{\text{if } f=0} \begin{cases} |X_I| = \frac{a_e \cdot |X|}{\sqrt{X^2+Y^2}} \\ |Y_I| = \frac{b_e \cdot |Y|}{\sqrt{X^2+Y^2}} \end{cases} \tag{15}$$

This new formulation shows that when $f = 0$ (i.e. $a_e = b_e$), the orthogonal contacting point exists, and hence, does not impose any numerical singularities for circles. Moreover, if equation 15 at $f = 0$ and radius, $R = a_e = b_e$ is employed, equation 13 simplifies to: $D_h(X,Y) = \sqrt{X^2+Y^2} - R$, which is the true geometric distance of a point to a circle. Another interesting observation from equation 15 is that when $X = 0$ (points on the minor axis), or when $Y = 0$ and $|X| \geq f$ (points on the major axis beyond the focal points), even though a confocal hyperbola does not technically exist, the equations of the intersecting point are still deterministic, and the distance functions reduce to the true geometric distance through the following simplified formulas:

$$\begin{cases} D_{hX}(X,Y) = 0 \qquad D_{hY}(X,Y) = |Y| - b_e \qquad \text{for } X = 0 \\ D_{hX}(X,Y) = |X| - a_e \qquad D_{hY}(X,Y) = 0 \qquad \text{for } Y = 0 \text{ and } |X| \geq f \end{cases} \tag{16}$$

Considering the presented formulations, equations 13, 15 and 16 can now be employed to find the geometric parameter vector $\rho = (x_c, y_c, a_e, b_e, \theta)$ by minimizing:

$$\min_{\rho} \sum_{i=1}^{N} \left(D_h(x_i, y_i)\right)^2 = \sum_{i=1}^{N} (D_{hX}(x_i, y_i)^2 + D_{hY}(x_i, y_i)^2) \quad (17)$$

where equation 17 is a non-linear problem, which must be solved using an iterative heuristic method such as Levenberg- Marquardt. To this end, **Algorithm 1: Ellipse Fitting using Confocal Hyperbola** is developed as follows:

1- Find an initial estimate of the geometric parameters $\rho = (x_c, y_c, a_e, b_e, \theta)^T$ of the ellipse using one of the methods presented in Section 2. Here, we will use Halir's numerically stable method, which is fast and guarantees an ellipse, as the initial parameter approximation.

2- Perform the Levenberg-Marquardt algorithm with maximum number of iterations, $N_T$, and damping parameter, $\lambda$, its decrement, $\gamma$, and increment, $v = v_0$, as follows:

a) For each point, calculate the distance function, $D_h$ (equations 15 and 16), and Jacobian, $JD_h$ (equations A6 and A7 of Appendix I), along with the sum of the squared distances, $SD = \sum_{i=1}^{N} {D_h}_i^{\,2}$, using the geometric parameters, $\rho$;

b) Estimate the new geometric parameters using the following equation:

$$\rho_N = \rho - (JD_h^T JD_h + \lambda . I_{5\times5})^{-1} JD_h^T D_h \quad (18)$$

c) Calculate the sum of the squared distances, $SD_N$, using the new geometric parameters, $\rho_N$:

- If $SD_N = SD$: retain the solution as the optimal and exit the algorithm; else
- If $SD_N > SD$: $\lambda = v\lambda$ and $v = v^2$; else
- If $SD_N < SD$: $\lambda = \frac{\lambda}{\gamma}$, $v = v_0$ and $\rho = \rho_N$;

d) If the number of iterations equals $N_T + 1$ exit the algorithm and retain the solution with the lowest $SD$; else return to step 2-a.

It is worth mentioning that there are many variations of the Levenberg- Marquardt algorithm, some include strategies such as parameter-wise re-weighting of the damping parameter $\lambda$ [30] as well as geodetic acceleration [31], which may also be used to solve equation 17. The proposed implementation of the Levenberg-Marquardt algorithm was, however, sufficient to achieve the results presented in Sections 5 and 6 in all cases. Another important note is that equation 17 can also be used to derive the best fit hyperbola to a given set of points using Algorithm 1 subject to the parameter vector $\rho_h = (x_c, y_c, a_h, b_h, \theta)$. Fitting a hyperbola to a set of points will, however, not be considered in this manuscript. The closed formulation for the Jacobian, $JD_h$, of the distance function is provided in Appendix I.

*3.1. Time Complexity Analysis of Algorithm 1*

Given $N$ input data points, the time complexity of Algorithm 1 can be determined from the following steps:

1- Halir's initial ellipse parameter estimation is solved in $\sim O(N)$ time.

2- For each data point, the confocal hyperbola distance, $D_h$, requires 14 multiplications, 2 divisions, 7 subtractions, 6 summations, 4 square roots, and two absolute values, where each operation has constant time complexity $\sim O(1)$. Therefore, the complexity for all data points is $\sim O(N)$.

3- For each data point, the Jacobian matrix, $JD_h$, requires 60 multiplications, 15 divisions, 14 subtractions, 12 summations, 4 square roots, and two absolute values, where each operation has constant time complexity $\sim O(1)$. Therefore, the complexity for all data points is $\sim O(N)$.

4- $\rho_N$ consists of the following computations is eventually in the order of $\sim O(N)$:

   a) $P_{5\times5} = JD_{h\,5\times N}^{T} JD_{h_{N\times5}}$ is solved in $O(25N)$;

   b) $Q_{5\times5} = P_{5\times5} + \lambda . I_{5\times5}$ has 25 summations $O(25)$;

   c) The inverse ${Q_{5\times5}}^{-1}$ is $O(5^3)$;

   d) $S_{5\times1} = JD_{h\,5\times N}^{T}\, D_{h_{N\times1}}$ is $O(5N)$;

   e) The multiplication $L_{5\times1} = {Q_{5\times5}}^{-1} . S_{5\times1}$ is $O(25)$;

   f) $\rho_N = \rho - L_{5\times1}$ is five subtractions $O(5)$.

Therefore, every iteration can be solved in linear $\sim O(N)$ time. Given the maximum number of iterations, $N_T$, the maximum time complexity will be in the order of $\sim O(N_T . N)$. Empirically, the expected number of iterations to convergence was shown to be around 9 iterations (see Figure 10b). Furthermore, the computation time (not complexity) of the confocal hyperbola distance, $D_h$, as well as the Jacobian matrix, $JD_h$, i.e. steps 2 and 3 above, can be significantly reduced, if computer code vectorization [32] is used instead of a "for loop".

## 4. Experimental Design

In this manuscript, two main categories of experiments are designed, namely, simulation, and real-world. Simulation-based experiments were designed to validate the results of the ellipse fitting compared to other state-of-the-art and established methods in data resembling edge points in images. The real-world experiments were designed to evaluate the ability of different ellipse fitting methods in recovering the geometric parameters of cylinders in 3-dimensional point clouds. The two classes are explained in more detail in the following. All computations were carried out using a computer with AMD Ryzen 5-2600X CPU, 64GB RAM, and 1TB SSD NVME storage. The summary description of the designed experiments is provided in Table 1.

**Table 1:** Summary of the designed experiments

| Experiment Description | | Purpose |
|---|---|---|
| **Simulation-based experiments** | **Point to ellipse distance measure** | Comparing the confocal hyperbola distance to Algebraic, Sampson , Harker [24], and Geometric [11] |
| | **Overall ellipse fitting evaluation** | Comparing overall RMSE and P-Error achieved by Algorithm 1 to the methods of Halir [16], Taubin [19], Kanatani [20], Szpak [12] and Ahn [27] |
| | **Parameter specific ellipse fitting evaluation** | Comparing the P-Error achieved by Algorithm 1 to the methods of Halir [16], Taubin [19], Kanatani [20], Szpak [12] and Ahn [27] as a function of:<br>(i) rotation angle;<br>(ii) aspect ratio;<br>(iii) noise; and<br>(iv) spanning arc angle |
| **Real-world experiments** | | Evaluating the ability of ellipse fitting using algorithm 1, compared to the methods of Halir [16], Szpak [12] and Ahn [27] in recovering geometric parameters of cylinders from point clouds |

### *4.1. Metrics for validation of results*

For the simulation-based experiments, two metrics are used: (i) the root mean squared error (RMSE) of the best fit ellipse; and (ii) the L2-norm of the estimated best fit parameter vector from the ground truth divided by the L2-norm of the ground truth parameter vector (abbreviated here as P-Error and reported in percentages). The RMSE of the best fit will be calculated using the method of Chernov [11], which is currently the most reliable method for projecting points onto ellipses, and consequentially determining the actual geometric distance from points to the ellipse. For the real-world experiments, the L2-norm of the estimated cylinder parameters from the ground truth parameters, namely center, radius, and axis angle, are used. The ground truth geometric parameters of the cylinders from point clouds are determined using the robust cylinder extraction method proposed by Maalek [8].

### *4.2. Simulation-based experiments*

In this section, the method in which elliptic points are generated along with the various simulated configurations used to evaluate different ellipse fitting methods are explained.

#### *4.2.1. Simulating elliptic points*

An important consideration for the simulated experiments is the method in which simulated ellipses are generated. Here, our focus for simulation experiments is producing elliptic points resembling edge points in images (three-dimensional (3D) point clouds will be considered in the real-world experiments). The edge points in images contain two important characteristics. First, the larger the ellipse (or any shape for that matter), the higher the number of ellipse edge points will be. In other words, an un-occluded circle with radius of 3 pixels cannot comprise of say 100-pixel points. Second, the discrete pixelated nature of digital images ensures that an extracted edge point will always fall on an image pixel. In other words, even if no random measurement error exists, the edge points representing an ellipse in an image will fall on the closest pixel and not necessarily exactly on the ellipse. The effects of digitization along with those induced by the number of points with respect to the scale of an ellipse must be considered when simulating ellipses. A process is, hence, formulated to automatically determine the pixel points representing the ellipse in images using only the geometric parameters of the ellipse. To this end, **Algorithm 2: Simulating Elliptical Edge Points**, given the geometric parameter vector $\rho = (x_c, y_c, a, b, \theta)^T$, and the polar arc angles for the beginning, $\alpha_s$, and the end, $\alpha_F$, of the ellipse, subjected to random measurement error, $\mathcal{N}(0, \sigma^2)$, is developed as follows:

1- Find the orthogonal circumscribed rectangle to the ellipse in digital coordinates using the following:

$$\begin{cases} x_{rm} = \lfloor x_c - \sqrt{a^2 \cos\theta + b^2 \sin\theta} \rfloor, & x_{rM} = \lceil x_c + \sqrt{a^2 \cos\theta + b^2 \sin\theta} \rceil \\ y_{rm} = \lfloor y_c - \sqrt{b^2 \cos\theta + a^2 \sin\theta} \rfloor, & y_{rM} = \lceil y_c + \sqrt{b^2 \cos\theta + a^2 \sin\theta} \rceil \end{cases} \quad (19)$$

Where $(x_{rm}, x_{rM})$ and $(y_{rm}, y_{rM})$ are the minimum and maximum of the limits of the orthogonal circumscribed rectangle in $x$ and $y$, respectively, and $\lfloor . \rfloor$ and $\lceil . \rceil$ denote the floor and ceiling functions, respectively.

2- Generate a grid of $k = x_{rM} - x_{rm} + 1$ by $l = y_{rM} - y_{rm} + 1$ points between the rectangular limits of $(x_{rm}, x_{rM})$ and $(y_{rm}, y_{rM})$ with one pixel spacing in the $x$ and $y$ direction, respectively.

3- Determine the orthogonal contacting point $(x_{PI}, y_{PI})$ of each pixel point $(x_P, y_P)$ within the generated grid from the ellipse using the method of Chernov [11].

4- Find the pixels satisfying the following conditions:

$$\begin{cases} |x_{PI} - x_P| \leq 0.5 \\ |y_{PI} - y_P| \leq 0.5 \end{cases} \quad (20)$$

5- Find the orthogonal contacting points $(x_{PI}, y_{PI})$ of step 4 whose polar angles are between $\alpha_s$ and $\alpha_F$.

6- Subject the $(x_{PI}, y_{PI})$ coordinates satisfying steps 4 and 5 to additive random measurement error $\mathcal{N}(0, \sigma^2)$. $\sigma$ is in pixel coordinates, and $x_{PI}$ and $y_{PI}$ are considered independent and identically distributed (i.i.d).

7- Find the closest pixel grid to the coordinate points of step 6.

Figure 2a and 2b illustrate the simulated pixels of an ellipse with the geometric parameter vector, $(1,3,15,10,30^\circ)^T$, with no measurement error, and subjected to random measurement error, $\mathcal{N}(0,1)$, respectively.

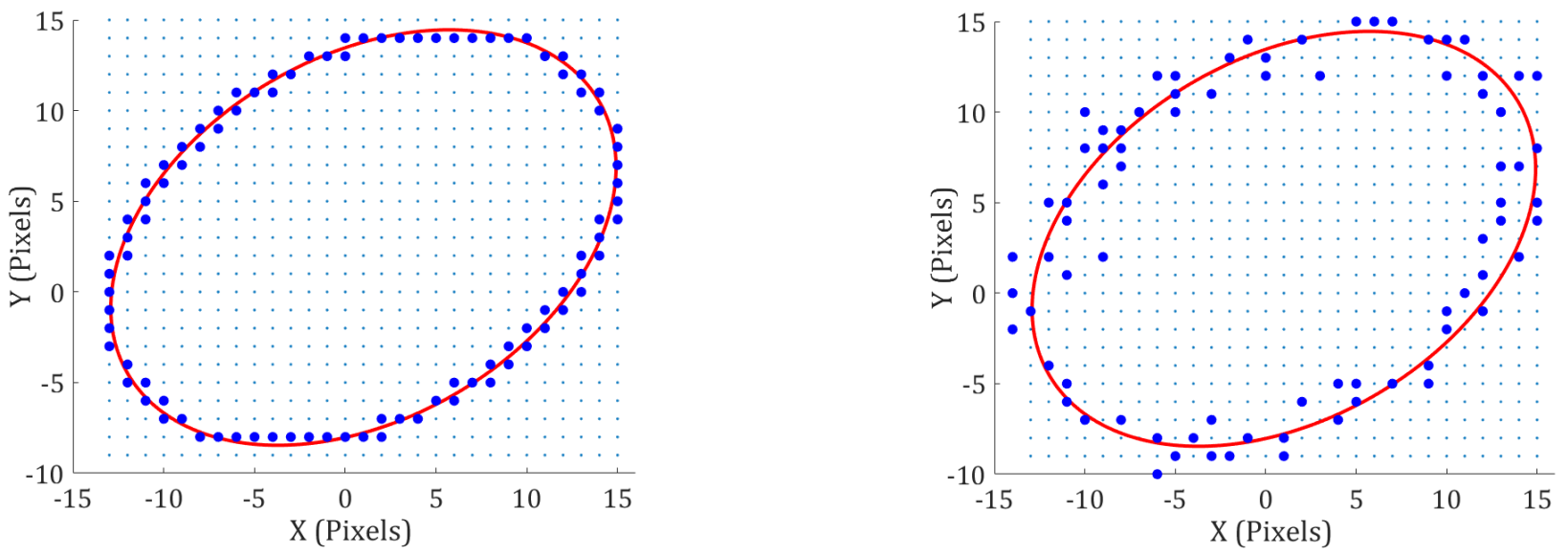

**(a)** **(b)**

**Figure 2:** Simulated ellipse with parameter vector $(1,3,15,10,30^{\circ})^T$: a) ellipse with no random measurement errors (only the pixel rounding effect); b) ellipse with random measurement error $\mathcal{N}(0,1)$;

*4.2.2. Ellipse simulation configurations*

Predicated on the proposed methods, two main categories of simulation-based experiments are carried out, namely, point to ellipse distance measure, and ellipse fitting evaluation (each with two sub-categories), as follows:

1- **Point to ellipse distance measure:** is designed to numerically identify the best method amongst various point to ellipse distance measures to accurately predict the ground truth geometric distance, obtained by Chernov [11]. To this end, 10,000 ellipses were randomly generated using Algorithm 2, subject to the following conditions:

   a) **rotation angle,** randomly chosen between 0 and $\pi$;

   b) **aspect ratio,** randomly chosen between 1.5 and 4;

   c) **measurement standard deviation,** randomly chosen between 0.5 and 15 pixels;

   d) **semi-minor,** was chosen as 50 pixels; and

   e) **ellipse spanning arc ($\boldsymbol{\alpha_F - \alpha_S}$),** was chosen as $\frac{\pi}{2}$ (i.e. quarter ellipse).

   The ellipse spanning arc was selected as a quarter of an ellipse (rather than say a full ellipse) since the point to ellipse measurements are not expected to drastically change between different quadrants of the same simulated ellipse. For each measurement, the distance of point to ellipse is calculated using the confocal hyperbola, algebraic, Sampson (see Section 2), and Harker and O'Leary [24] . The absolute deviation of the distances for each method from that obtained by Chernov is reported.

2- **Overall ellipse fitting evaluation:** this experiment is designed to evaluate the effectiveness of the proposed ellipse fitting method using confocal hyperbola (Section 3) compared to the ellipse fitting methods of Ahn [27], Szpak [12], Kanatani [23], Taubin [19], and Halir [16]. Similar to experiment 1, 10,000 randomly generated ellipses were used to determine which ellipse fitting method produces the fit with the lowest RMSE along with P-Error (see Section 4.1). The ellipse parameter configurations were as follows:

a) **rotation angle,** randomly chosen between 0 and $\pi$;

b) **aspect ratio,** randomly chosen between 1.5 and 4;

c) **measurement standard deviation,** randomly chosen between 0.5 and 5 pixels;

d) **semi-minor,** randomly chosen between 10 and 100 pixels; and

e) **ellipse spanning arc,** randomly chosen from either full ($2\pi$), $\frac{3\pi}{2}$, or half ($\pi$).

For each randomly-selected configuration, 500 different simulations were carried out to capture the behavior of random measurement errors. The results of the mean, median and 95$^{th}$ percentile for the 500 simulations in each configuration, achieved by each method, were reported. In addition, the methods of Ahn, Szpak and ours are iterative, whereas the rest are direct. To provide a fair comparison between the different iterative methods, similar iteration-termination criteria are used. For all methods, Halir's ellipse fitting is used to estimate the initial ellipse parameters. The damping parameter (step size for Ahn's method), damping increment, damping decrement, and maximum number of iterations were set as $\lambda = 0.5$, $v = 10$, $\gamma = 3$, and 50, respectively.

3- **Parameter specific ellipse fitting evaluation:** The last set of simulated experiments involves the assessment of the impact of variables, namely, rotation angle, aspect ratio, measurement standard deviation, and spanning elliptic arc, on the accuracy of the estimated parameters of the best fit ellipse (P-Error) using the methods of Ahn [27], Szpak [12], Kanatani [23], Taubin [19], and Halir [16], and Prasad [33]. To this end, a base parameter setup is used, and only one parameter is changed at a time with the configurations, shown in Table 2. Similar to the previous experiment, to correctly capture the behavior of random measurement errors on the performance of each ellipse fitting method, each parameter configuration was simulated 500 times and the mean of the P-Error was recorded.

**Table 2:** Configurations for the parameter specific ellipse fitting evaluation

| Parameter Type | Base Configuration | Parameter Change | | |
|---|---|---|---|---|
| | | From | To | Increments |
| **Rotation angle** | $\frac{\pi}{4}$ | 0 | $\pi$ | $\frac{\pi}{20}$ |
| **Aspect ratio** | 2 | 1 | 4 | $\frac{3}{20}$ |
| **Noise** | 2 | 0 | 5 | 0.25 |
| **Spanning arc angle** | $2\pi$ | $\frac{\pi}{2}$ | $2\pi$ | $\frac{3\pi}{40}$ |

### 4.3. Real-world experiments

In the early 1800, Dandelin showed that the intersection between a plane and a cylinder produces an ellipse. The geometric parameters of the cylinder can be retrieved as a function of the ellipse parameters as follows [34]:

$$\begin{cases} R = b_e \\ \delta = \cos^{-1} \frac{b_e}{a_e} \\ C_e = Rot.\left(C_c + \frac{\langle P_0 - C_c, \vec{n} \rangle}{\langle \vec{a}, \vec{n} \rangle} \vec{a}\right) \end{cases} \tag{21}$$

where $R$, $C_c$ and $\vec{a}$ are the radius, a point on the axis (center), and vector of cylinder's axis, respectively, $\vec{n}$ and $P_0 = (x_0, y_0, z_0)$ are the normal vector to and a point on the intersecting plane, respectively, $C_e = (x_e, y_e, z_0)^T$, $a_e$ and $b_e$ are the center, semi-major length and semi-minor length of the ellipse in the coordinate system of the intersecting plane, respectively, $Rot$ is the rotation matrix that takes the plane's normal vector to vector $(0,0,1)^T$, and $\delta$ is the angle between the plane's normal vector and the cylinder's axis. From this formulation, it is possible to approximate the parameters of a cylinder using the parameters of an ellipse created through an intersecting plane. To this end, this experiment is designed to evaluate the ability of different ellipse fitting methods in recovering the geometric parameters in three separate real-world cylindrical point clouds. The selected point clouds are half cylinders of mechanical pipes, acquired using a Leica HDS6100 terrestrial laser scanner (Figures 3a and 3b), as well as 3D reconstructed point cloud using twenty iPhone 11 4K images (Figure 3c). The point clouds represent three noise levels and point densities with a normalized RMSE of approximately 0.010, 0.015 and 0.020 (for the point clouds of Figures 3a, 3b and 3c, respectively). The normalized RMSE is defined, here, as the RMSE divided by the square root of the area of the circle (or ellipse). The outlier removal, and the ground truth cylinder parameter estimation were carried out using the robust cylinder extraction and fitting proposed in [8]. For each point cloud, 10,000 random planes passing through a randomly selected point on the cylinder were generated (Figure 3d-top). The closest 50 points on the cylinder to each plane -within 1mm distance from the plane- were identified (Figure 3d-bottom). The points are then rotated such that the normal to the plane is parallel to $(0,0,1)^T$. The best fit ellipse using the methods of Halir [16], Szpak [12], Ahn [27], and ours is applied on the $x - y$ coordinates of the rotated points (only methods guaranteeing an ellipse were considered). From the best-fit ellipse parameters estimated by each method, the geometric parameters of the original cylinder were computed using equation 21, and the errors between the estimated parameters and the ground truth were calculated. Here for clarity, Figure 3d-bottom shows the ground truth ellipse in randomly generated plane (Figure 3d-top), computed from back calculating equation 21 using the ground truth geometric cylinder parameters.

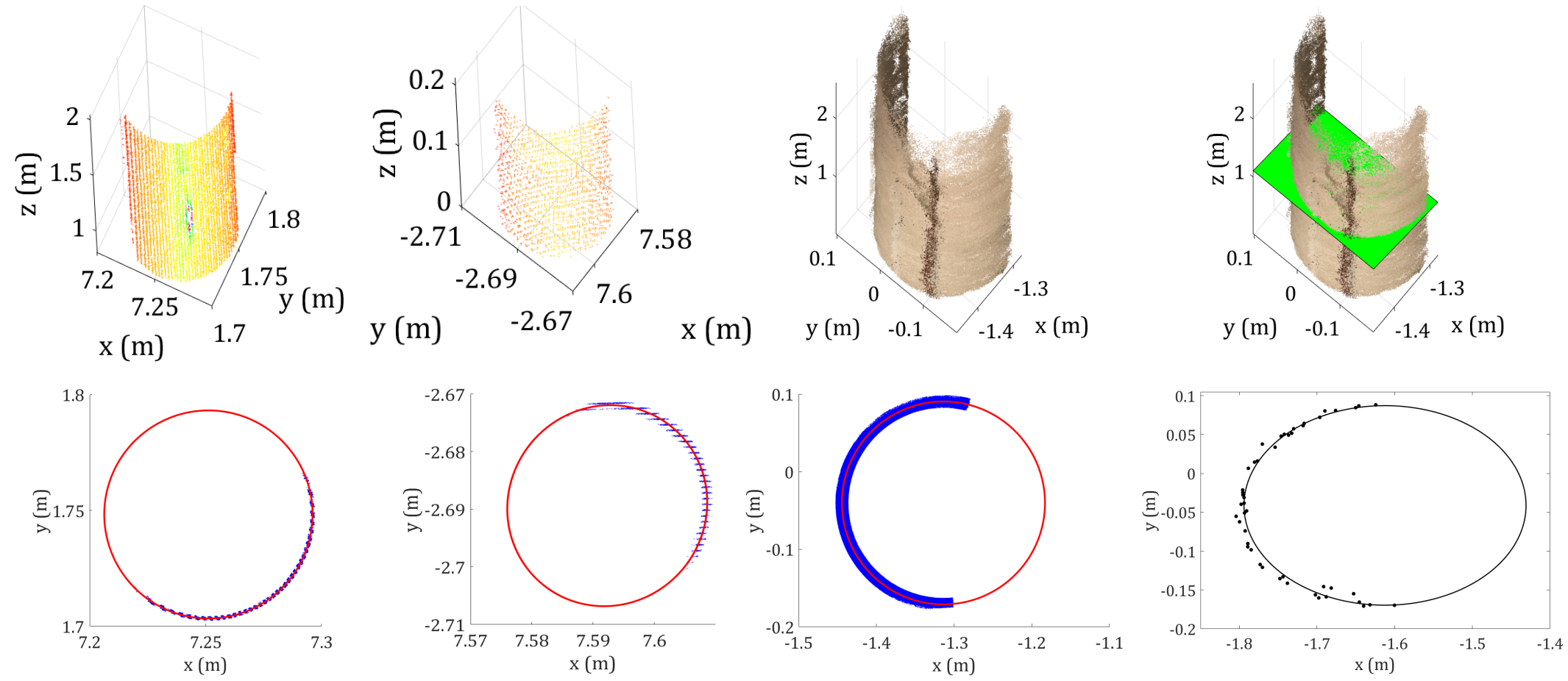

**(a)** **(b)** **(c)** **(d)**

**Figure 3:** Sample point clouds of cylindrical pipes extracted using [8], oriented such that the cylinders' axis is parallel to $(0,0,1)^T$: a) TLS with low noise and medium point density; b) TLS with medium noise and low point density; and c) 3D reconstruction of images with high noise and high point density. d) Example of a generated plane with random orientation (green-top) passing through the cylindrical point cloud of Figure 3c, and the consequential elliptical points in the vicinity of the intersecting plane (bottom).

## 5. Simulated experiment results

### *5.1. Point to ellipse distance measure*

In Section 3, a new ellipse fitting method using the confocal hyperbola distance was introduced. Here, we schematically and empirically show the agreement of the confocal hyperbola distance to the ground truth geometric distance, acquired using Chernov's method [11]. The behavior of the confocal hyperbola distance function is also compared to the algebraic, Sampson, and Harker and O'Leary distance approximations.

#### *5.1.1. Visual behavior of distance approximations*

Before the numerical results are presented, it is important to also study the behavior of the point to ellipse distance measurements using the confocal hyperbola, algebraic, Sampson, Harker and O'Leary [24], and Chernov [11] (the ground truth), particularly in ill-positioned points relative to the ellipse (i.e. points on the major or minor axes). To this end, consider an ellipse with geometric parameter vector of $(0,0,5,3,0)^T$. Figure 4 shows the behavior of the distance functions of points in three different arrangements, namely, points on the major axis (Figure 4a: $\{x \in [-10,10], y = 0\}$), points on the minor axis (Figure 4b: $\{y \in [-10,10], x = 0\}$), and points on the line $x = y$ (Figure 4c: $\{x, y \in [-10,10]\}$). Figures 4a, 4b and 4c all show that the Sampson distance tends to infinity as the point coordinates approach the center of the ellipse. Figure 4a shows that the Harker and O'Leary's distance also tends to infinity as the point coordinate approaches the focal points of the ellipse (in our case points $x = \pm 4$ and $y = 0$).

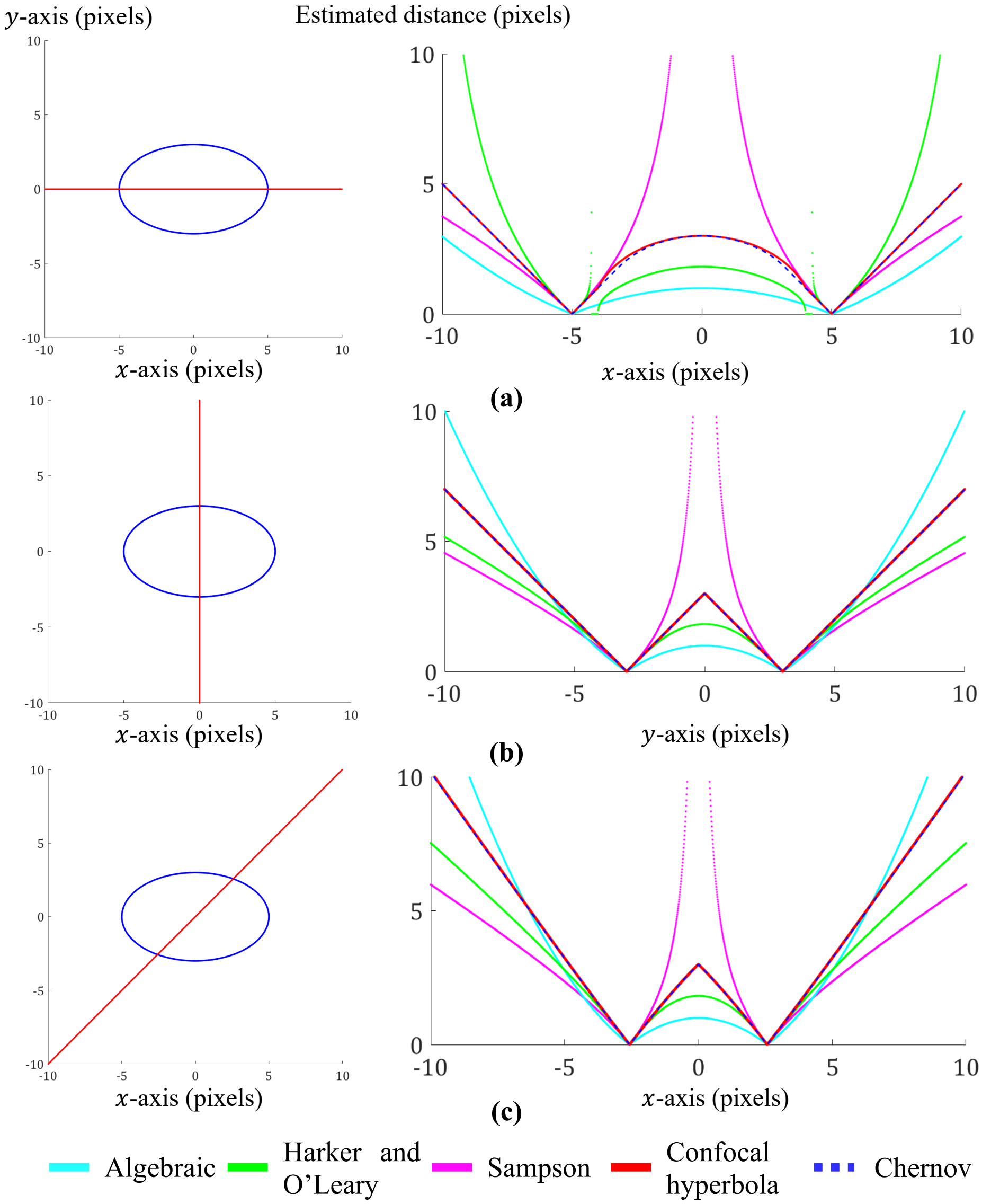

**Figure 4:** Behavior of the distance functions of algebraic (cyan), Sampson (magenta), Harker and O'Leary (green), confocal hyperbola (red), and Chernov (blue) presented on the right for the distance of the red points to the ellipse (with geometric parameters $(0,0,5,3,0)$) for the arrangements shown on the left: a) points on the major axis; b) points on the minor axis; and c) points on the line $x = y$.

Another interesting observation from Figure 4 is that none of the distances obtained using algebraic, Harker and O'Leary, and Sampson, correctly predicts the pattern of the true geometric distance (i.e. Chernov's method shown in blue). In fact, the distances of Harker and O'Leary, Sampson, and algebraic appear to move even farther away from the actual geometric distance as the point moves farther from the ellipse. Especially for the case of Sampson-based methods, the visual observations corroborate that Sampson distance may only be appropriate for data with smaller noise levels as was mentioned in [12]. Therefore, careful attention must be given when integrating the Sampson distance into ellipse fitting cost functions, especially at larger noise levels (see Figure 3 of [35] as an example). The confocal hyperbola method on the other hand, is almost indistinguishable from the true geometric distance in the presented three arrangements (the confocal hyperbola, shown in red lies almost exactly underneath Chernov's distance shown in blue dashed line). From Figure 4, it can also be observed that the confocal hyperbola, Chernov and the algebraic distances do not impose any singularities at these ill-positioned point arrangements, unlike the distance functions of Sampson, and Harker and O'Leary.

*5.1.2. Numerical comparison*

As observed, the visual assessment of the distance functions shows the advantages of using the confocal hyperbola method compared to other methods in predicting the behavior of the geometric distance. The results of the simulation presented in the following are, hence, used to numerically quantify the extent of this advantageous performance. To this end, 10,000 different ellipse configurations were simulated as per the properties presented in Section 4.2.2. For each measurement, the absolute difference of the point to ellipse distances using the algebraic, Sampson, Harker and O'Leary and confocal hyperbola, from that of the ground truth geometric distance was calculated. Figure 5a shows the results of the simulated measurements. The horizontal axis represents the absolute difference of each method from the ground truth in pixels, and the vertical axis represents the cumulative probability. It can be clearly observed that the confocal hyperbola achieved significantly better results compared to the remaining methods. To provide some numerical perspective, the mean, median, and 95$^{th}$ percentile of the absolute distances obtained by each method from the ground truth are presented in Table 3. It can be observed that the 95$^{th}$ percentile of all measurements using the confocal hyperbola fall within $\pm 0.05$ pixels from the ground truth. The confocal hyperbola is, in fact, approximately 34, 73, and 358 times more accurate than the distance approximations of Harker and O'Leary, Sampson, and algebraic (on average). It can also be observed that 50% of the distances using the confocal hyperbola are less than $10^{-4}$ pixels from the ground truth, which is 99, 1136, and 26792 times better than those achieved by the best 50% of Harker and O'Leary, Sampson, and algebraic distances. In fact, the Sampson distance, which has been known as an excellent approximation to the geometric distance [12], is inferior to even Harker and O'Leary, which is also considerably inferior to the confocal hyperbola method.

**Table 3:** Mean, median and 95th percentile of the absolute difference of each method from the ground truth (pixels)

| **Method** | **Algebraic** | **Sampson** | **Harker and O'Leary** | **Confocal Hyperbola** |
|---|---|---|---|---|
| Mean | 5.81 | 0.77 | 0.31 | 0.01 |
| Median | 3.91 | 0.16 | 0.01 | 0.00 |
| 95th Percentile | 17.98 | 3.66 | 1.72 | 0.05 |

Figure 5b shows the computation times for each method w.r.t the number of points in the dataset (from 50 to 10,000 points). It was observed that the ground truth method of Chernov was on average approximately 900, 150, 4700, and 12700 times slower than confocal hyperbola, Harker and O'Leary, Sampson, and algebraic methods, respectively. Even though the confocal hyperbola achieved about 5 times slower computation compared to Sampson, it still computed the distance of 10,000 points in around 0.7 milliseconds, while providing highly accurate geometric distance approximations.

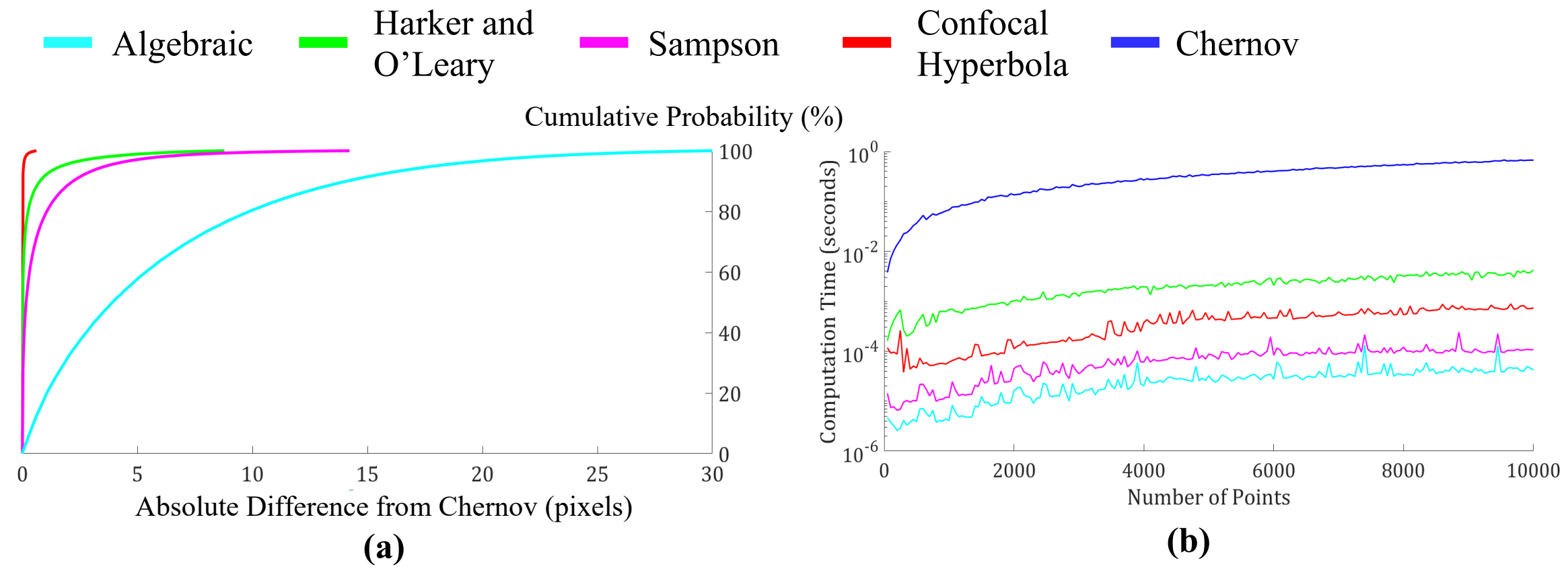

**Figure 5:** a) Absolute difference of the estimated distances using algebraic (cyan), Sampson (magenta), Harker and O'Leary (green), and confocal hyperbola (red) from the Chernov's method for all simulated measurements; b) computation time (in logarithmic scale) vs. number of points for different methods

### *5.2. Ellipse fitting evaluation*

In this section, first, the performance (RMSE and P-Error) of the proposed ellipse fitting using confocal hyperbola in comparison to the methods of Halir [16], Taubin [19], Kanatani [20], Szpak [12], and Ahn [27] was empirically evaluated. To this end, 10,000 random ellipses were simulated using the configurations presented in Section 4.2.2. The mean, median and 95th percentile of the RMSE and P-Error are provided in Table 4. The bolded numbers represent the best (lowest number) between all methods in each row. As shown, our ellipse fitting provides almost identical (and in some cases slightly better) RMSE and P-Error compared to the gold standard method of Ahn, and outperformed the remaining methods, particularly in the parameter estimation errors. Even though our method outperformed the rest, the difference in the RMSE between our method and the third best performing method, Szpak is still only in the order of $10^{-2}$ pixels, which may be negligible in practical settings. On the other hand, the P-Error, which represents the agreement between the estimated and original parameters, showed around 0.1%, 5.5%, 23.7%, 26.3%, and 107.3% relative improvement compared to Ahn, Szpak, Kanatani, Taubin and Halir, respectively (the average improvement amongst all measurements).

The results, shown in Table 4, provide an overall indication of the expected performance of our method compared to established state-of-the-art methods. It is, however, also important to assess the performance (particularly the P-Error) of all methods with respect to individual variables, such as, rotation angle, aspect ratio, noise level, and spanning arc. The results of these simulations are presented in the following subsections.

**Table 4:** Mean, median and 95<sup>th</sup> percentile of the RMSE and P-Error for 10,000 simulated ellipses

| **RMSE (pixels)** | | | | | | |
|---|---|---|---|---|---|---|
| **Methods** | **Halir** | **Taubin** | **Kanatani** | **Szpak** | **Ahn** | **Our Method** |
| Mean | 1.79 | 1.76 | 1.76 | 1.76 | **1.75** | **1.75** |
| Median | 1.78 | 1.77 | 1.77 | 1.76 | **1.76** | **1.76** |
| 95th Percentile | 2.92 | 2.86 | 2.86 | 2.86 | **2.85** | **2.85** |
| **P-Error (%)** | | | | | | |
| **Methods** | **Halir** | **Taubin** | **Kanatani** | **Szpak** | **Ahn** | **Our Method** |
| Mean | 1.048 | 0.665 | 0.637 | 0.608 | 0.570 | **0.569** |
| Median | 0.473 | 0.364 | 0.353 | 0.321 | **0.307** | **0.307** |
| 95th Percentile | 3.837 | 2.290 | 2.231 | 2.081 | 1.873 | **1.868** |

Note: The bold numbers represent the best method for a given row

*5.2.1. Impact of rotation*

Here, the impact of change in ellipse rotation angle from zero to $\pi$, in $\frac{\pi}{20}$ increments, on the P-Error for the methods of Halir, Taubin, Kanatani, Ahn, Szpak, Prasad, and confocal hyperbola was evaluated. The mean of the P-Error for all rotation angles using the different methods is provided in Table 5. As illustrated, our method achieved the best results followed very closely by Ahn's method. Figure 6 shows the results of the mean of the P-Error for the 500 separate simulations per ellipse configuration vs. the rotation angle. It can be observed that all methods, except for Prasad, remain almost constant with respect to the rotation angle. In other words, the rotation angle does not appear to impact the results of the ellipse fitting for Halir, Taubin, Kanatani, Ahn, Szpak and our method. The method of Prasad, on the other hand, is highly impacted (approximately 60% difference between the lowest and the highest value), especially at rotation angle equal to $\frac{\pi}{2}$ (i.e. vertical ellipse), which demonstrates that the parameters estimated using Prasad's method are not invariant to affine transformation of the data. The reason for this lies in the heart of the objective function of Prasad (see equation 44 of [33]), where the objective function is multiplied by the $y$ axis coordinates to enhance numerical stability of their method. Due to the systematic trend imposed through affine transformation of the data, the method of Prasad will no longer be included in the remaining comparisons.

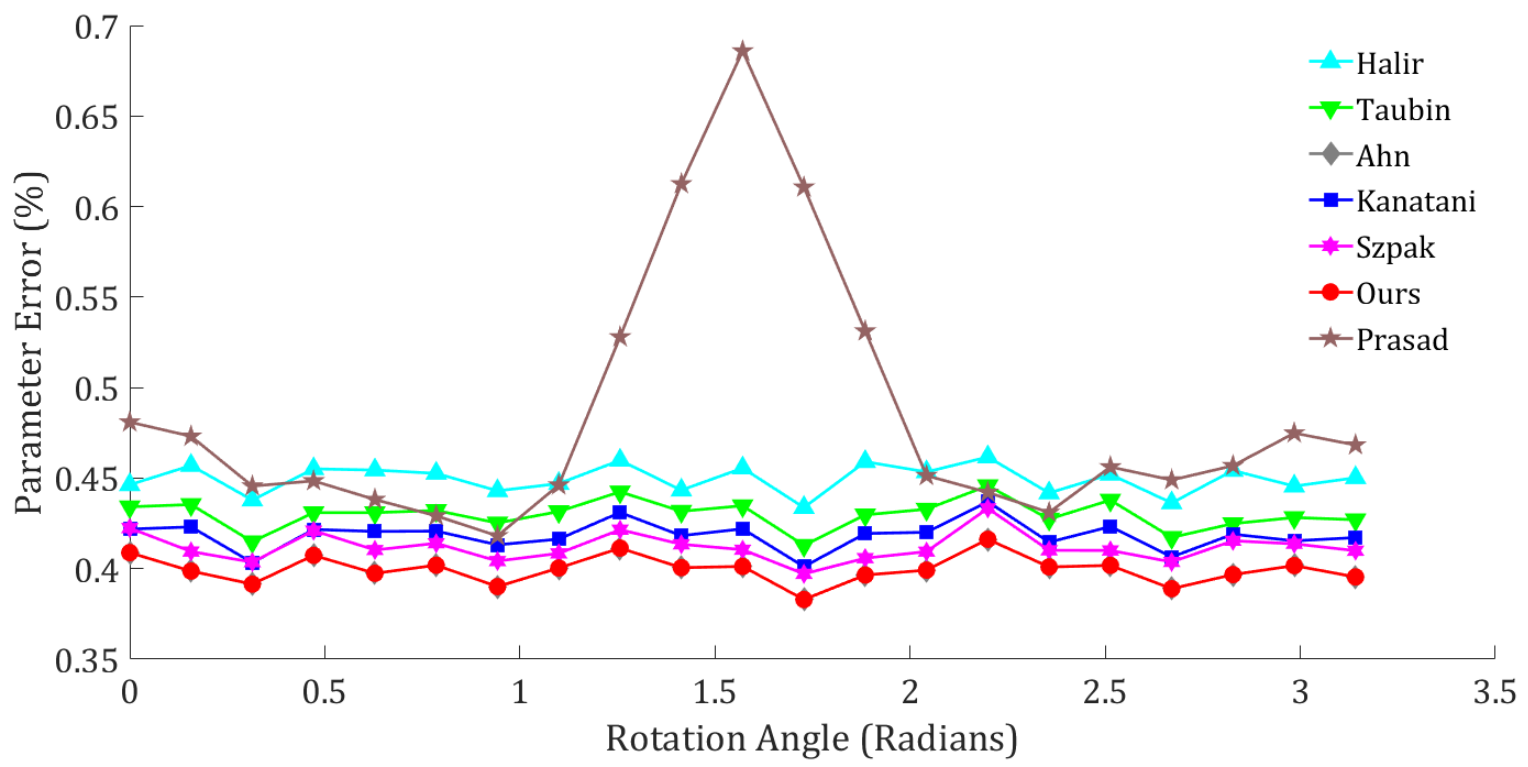

**Figure 6:** The mean of the parameter estimation errors for 500 simulation vs. the rotation angle (in radians)

*5.2.2. Impact of aspect ratio*

The aspect ratio was changed in 20 evenly spaced increments between 1 and 4. The mean of the P-Error for 500 simulations, obtained by the methods of Halir, Taubin, Kanatani, Ahn, Szpak, and confocal hyperbola are shown in Figure 7a. Three important observations were inferred from the provided chart. Firstly, all methods appear to behave similarly up until the aspect ratio of 1.5, where Halir's method starts to produce comparatively biased estimates of the best fit ellipses. Secondly, at around the aspect ratio of 2, the methods of Kanatani and Taubin start to produce comparatively biased results and appear impacted by the elongation of the ellipse. The methods, of Ahn, Szpak and ours appear to perform similarly with ours and Ahn's obtaining slightly better results on average compared to Szpak's (see Table 5). The last important observation is that all methods at the aspect ratio of 1 (i.e. a circle) appear to produce biased estimates of the geometric parameters. In fact, a relatively steep decline can be observed between the aspect ratio of 1 and 1.15. This demonstrates that the available ellipse fitting method may not necessarily be suitable for purely circle fitting applications. It is, hence, desirable to determine whether aspect ratios exist, where a circle fitting method can predict the model parameters more effectively. To this end, the ellipse fitting methods are compared to the reliable hyper-accurate circle fitting of Chernov [36] for aspect ratios between 1 and 1.012 in 0.001 increments. The results are shown in Figure 7b. It can be observed that for very small aspect ratios, it is possible to utilize the hyper accurate circle fitting method to predict the parameters of an ellipse. The circle fitting, however, starts to produce comparatively biased results at the aspect ratio of approximately above 1.005. This is since the semi-major and semi-minor of the simulated ellipses start to differ and a single radius estimation can no longer account for the differences between the semi-minor and semi-major lengths. Therefore, given the configurations used in this section, circle fitting might produce more reliable results for aspect ratio of below 1.005. In such cases, we can set $f = 0$, $\theta = 0$ and $a_e = b_e$ in our formulations (equations 13, 15, and A1-7) to treat equation 17 as a circle fitting problem by reducing the degrees of freedom from five parameters to three.

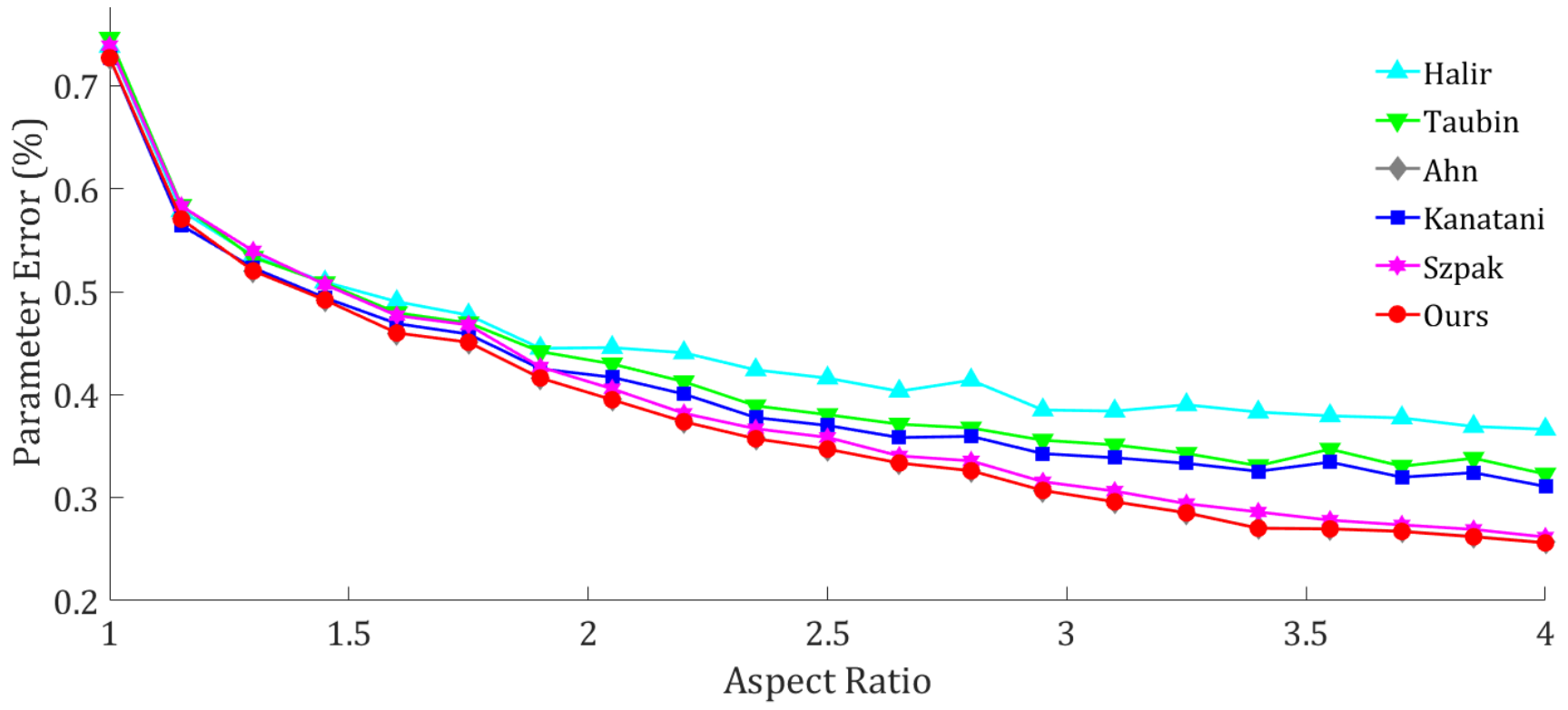

**(a)**

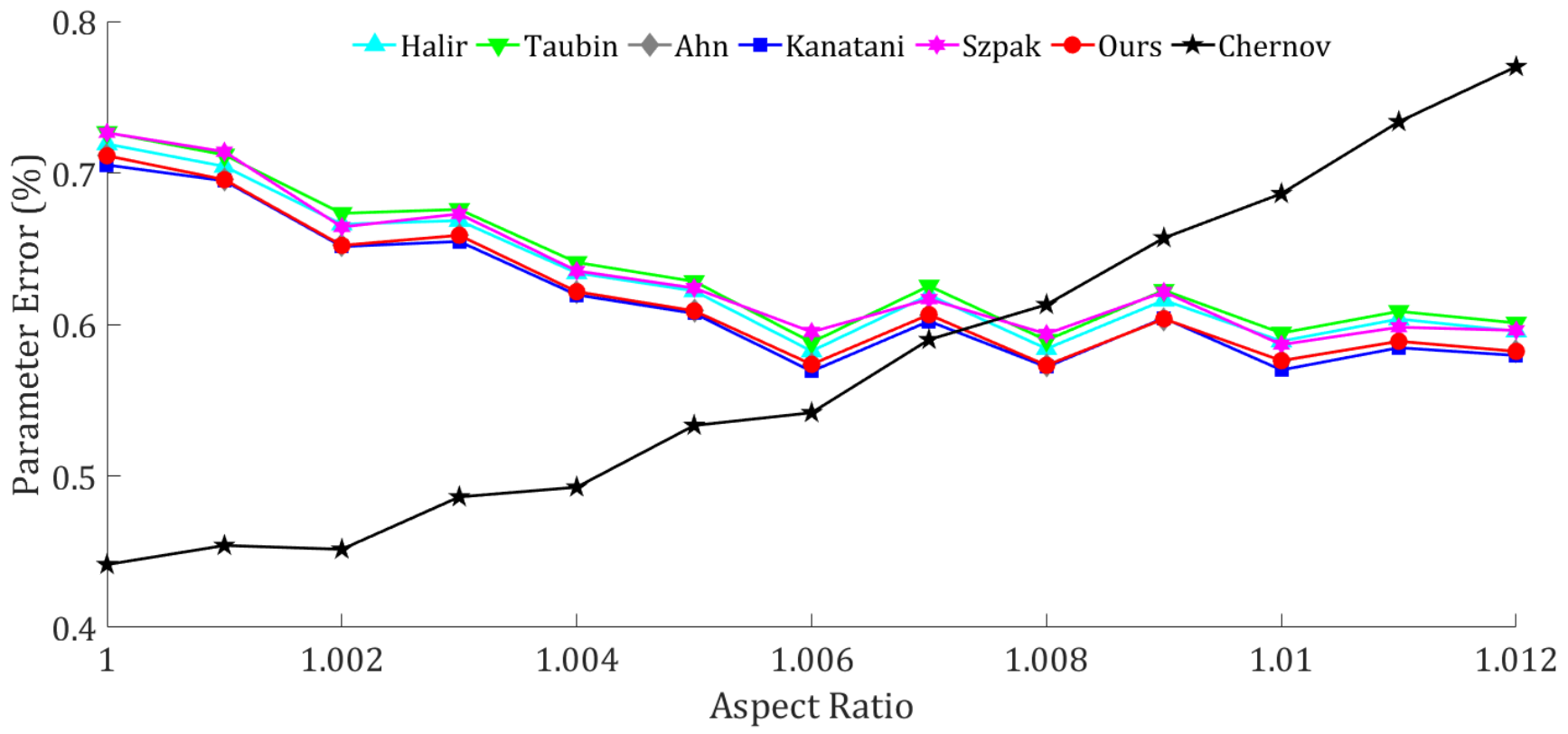

**(b)**

**Figure 7:** The mean of the parameter estimation errors for 500 simulations in each configuration vs. the aspect ratio: a) ellipse fitting methods; b) comparison of ellipse fitting with circle fitting method for small aspect ratios

*5.2.3. Impact of standard deviation (noise)*

Here, the measurement standard deviation is changed from 0 to 5 pixels at 0.25-pixel intervals. The results of the mean parameter error for each method vs. the measurement noise level is shown in Figure 8. Three main observations were drawn from the chart. Firstly, all methods appear to behave similarly up until around 1.5 pixels standard deviation, where Halir's method starts producing biased estimates of the parameters. The second observation is that at around 2.5 pixels standard deviation, the methods of Szpak and Taubin start producing relatively biased results compared to the methods of Kanatani, Ahn and ours. The last observation is that the methods of Ahn, Kanatani and ours are almost identical, and act very close to linear with respect to the standard deviation, with our method achieving an $R^2 = 0.99$ for the best fit line. The method of Kanatani was expected to behave favorably with respect to the standard deviation since the method was specifically designed to include bias corrections to the design and scatter matrices up to a second order error term. Ahn's method and ours appear to behave as favorably as that of Kanatani with respect to changes in measurement noise levels, while also guaranteeing an ellipse, which is one of the main disadvantages of Kanatani's method (see Section 2).

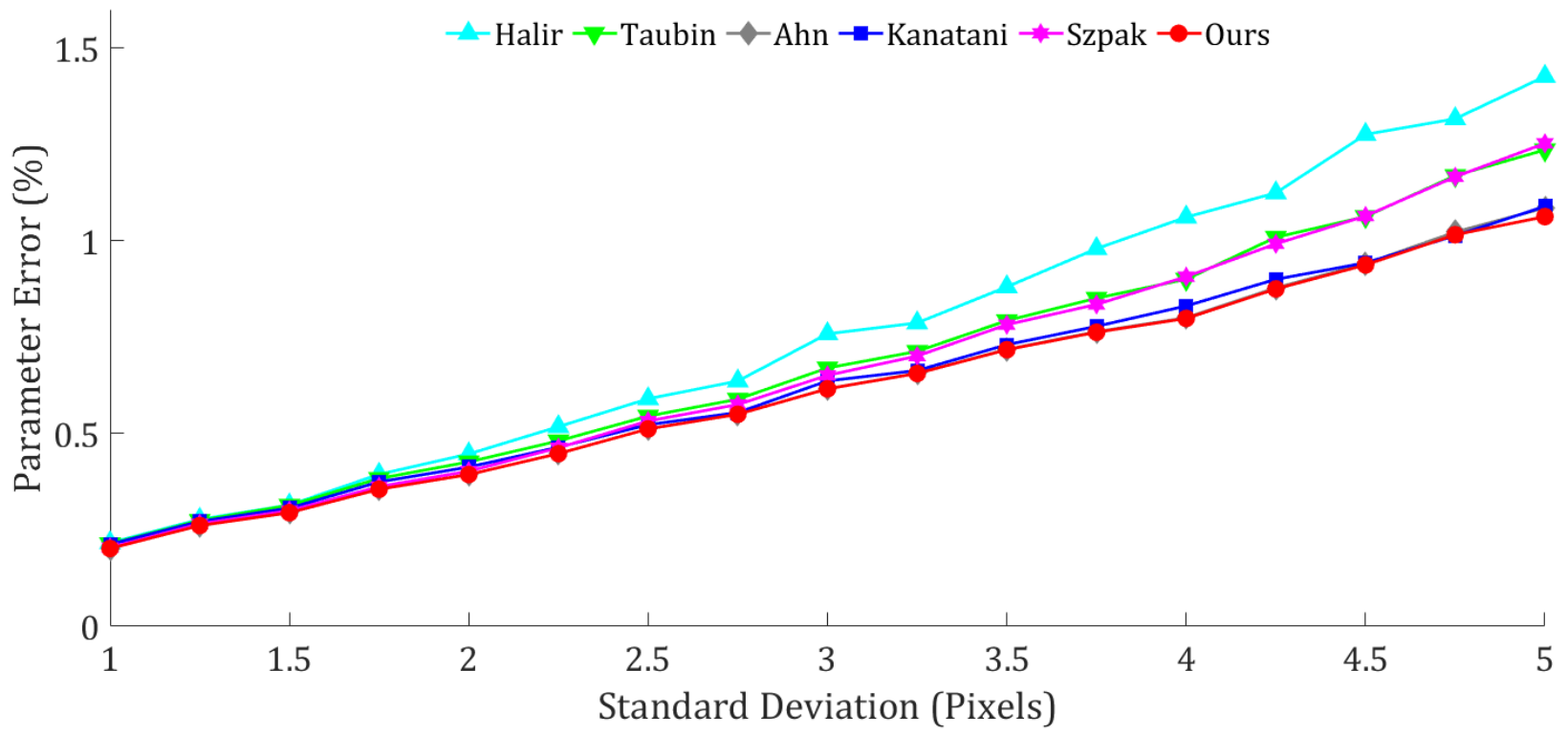

**Figure 8:** The mean of the parameter estimation errors for 500 simulations vs. the measurement standard deviation

*5.2.4. Impact of spanning arc angle*

The spanning arc angle is changed from $\frac{\pi}{2}$ to $2\pi$ in $\frac{3\pi}{40}$ intervals, and the results of the mean of the parameter estimation errors for each method is presented in Figure 9. It was observed that for spanning angles of over about $1.1\pi$ (~200°) all methods behaved similarly. At arc angles of approximately 200° and less, Halir's method started producing visibly biased results compared to the remaining five methods. The method of Kanatani and Taubin behave very similarly but start producing relatively biased results compared to ours, Szpak's and Ahn's at arc angles of around $0.65\pi$ (~120°). The method of Szpak and Ahn's performed very similarly at all arc angles with Ahn's method achieving a slightly better result on average in this category (Table 5). Our method started producing relatively better parameter estimation at arc angles of less than 120° compared to all other methods. Another interesting observation was that the relationship between the estimated parameter errors and arc angle closely followed a power function (e.g. $R^2 = 0.95$ for our method). In fact, using our method (the best performing method), the parameter estimation error for a quarter ellipse and a half ellipse were around 9.0% and 0.7% of the ground truth geometric parameters, respectively. Therefore, careful attention must be given when fitting ellipses to points on elliptic arcs of smaller than half an ellipse, particularly for applications pertaining to high precision metrology.

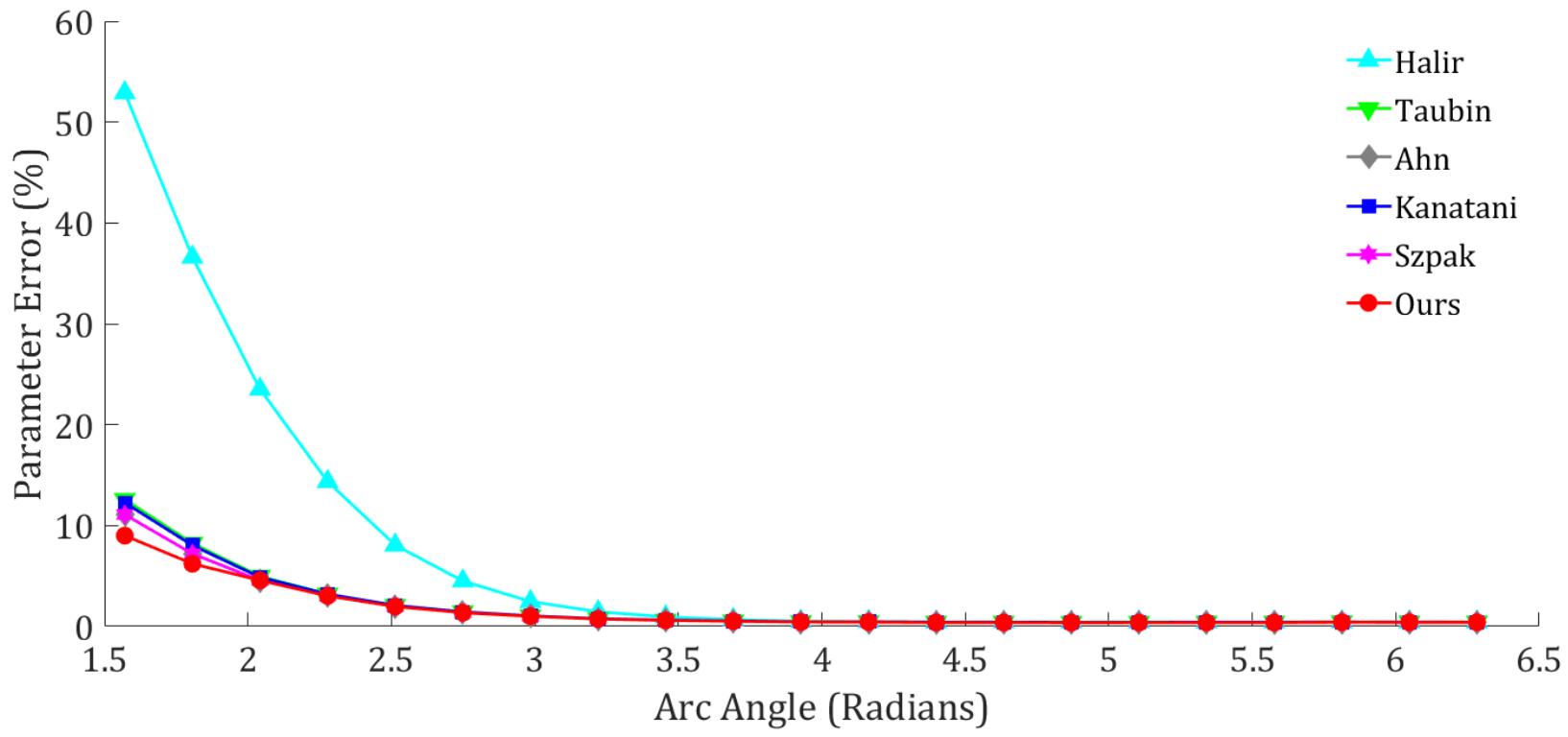

**Figure 9:** The mean of the parameter estimation errors for 500 simulations vs. the spanning arc polar angle (radians)

*5.2.5. Computation time and number of iterations*

In this section, the relative performance of our method is evaluated in terms of computation time along with the number of iterations to convergence. To this end, the base configuration of the ellipse simulation parameters of Section 4.2.2 was used. The number of points were increased from 200 to 4,000 in 200-point intervals, and at each interval 500 different simulations were performed. Figure 10a shows the mean of the computation time (for the 500 simulations at a given configuration) in logarithmic scale for each method vs. the number of points[1]. As observed, Ahn's, Szpak's and Kanatani's methods achieved the highest computation times, which were approximately 140-, 21- and 1.2- times slower than our method (average over all measurements). Our method on average was around 30 times slower than that of Halir's, which was also 1.4 times slower than Taubin's on average. To provide some perspective, ours, Szpak's and Ahn's took approximately 0.02, 0.66 and 3.91 seconds, respectively, to fitting an ellipse to 4,000 points. Not only did our method considerably outperform the other two iterative methods in terms of computation time, but it also outperformed Kanatani's direct method, particularly for number of points larger than 1,000.

The last comparison presented in this section is the number of iterations till convergence for the three iterative methods. Figure 10b shows the average number of iterations required for convergence vs. the number of points. As observed, Ahn's, Szpak's, and ours required on average approximately 19, 10 and 9 iterations, respectively to converge. An important observation is that Szpak's and ours converged after roughly the same number of iterations, yet our method achieved better computation times on average. This is since Szpak's method requires the evaluation of two objective functions at each iteration (see algorithm 2 of [12]). Not only does Ahn's method require more iterations (about twice as much as ours), which contributes to a higher computation time, but each iteration also takes longer to process. This is since in each iteration Ahn's method requires a separate generalized Gauss-Newton method to acquire the orthogonal contacting point to a given ellipse, a process which our method estimates very effectively through the confocal hyperbola distance function (see equation 13).

[1] The computation times shown in Figure 10a for iterative methods are based on the maximum iteration of 50, according to the convergence criteria provided in Section 4.2.2.

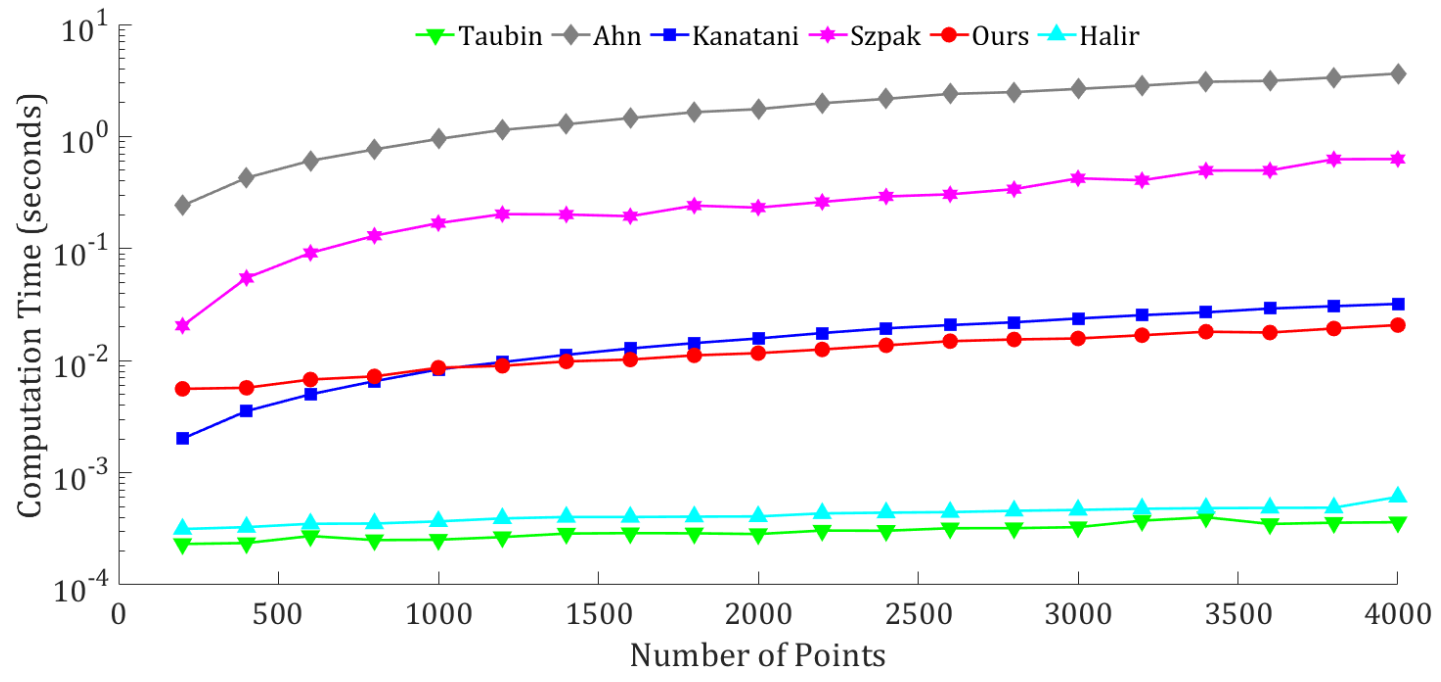

**(a)**

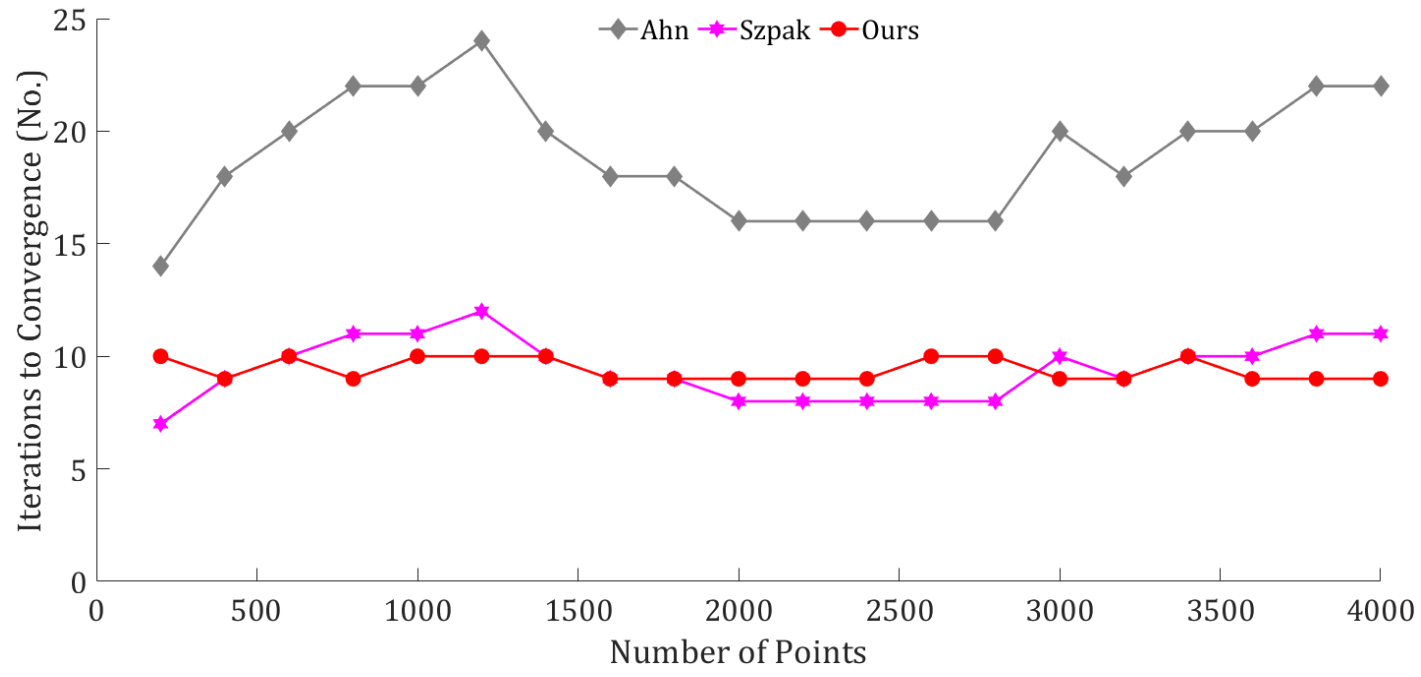

**(b)**

**Figure 10:** Simulation results for each method vs. number of points: a) computation time in seconds (logarithmic scale); and b) number of iterations to convergence

*5.2.6. Summary of best fit ellipse results*

The impact of four important variables, namely, rotation, aspect ratio, noise, and arc angle, on the parameter estimation using established ellipse fitting methods was evaluated. It was observed that the ellipse fitting methods of Halir, Taubin, Kanatani, Szpak, Ahn, and ours were not significantly impacted by the change on the rotation angle (see Figure 11a for reference in vertical ellipse). The methods of Szpak, Ahn and ours performed relatively more stable than the other methods with respect to the change in aspect ratio (see Figure 11b for behavior of the methods in aspect ratio of 4). It was also observed that the circle fitting of Chernov can provide more accurate parameter estimation in small aspect ratios (e.g. for the configurations used in the experiments less than 1.005). The methods of Ahn, Kanatani and ours performed relatively more stable than the other methods as a function of the measurement noise with a close to linear relationship between the estimated parameter errors and measurement noise (see Figure 11c for noise of 4 pixels). The methods of Ahn, Szpak and ours performed more stable than others, especially in arc angles of less than $\pi$, half ellipse (see Figure 11d for reference of half ellipses).

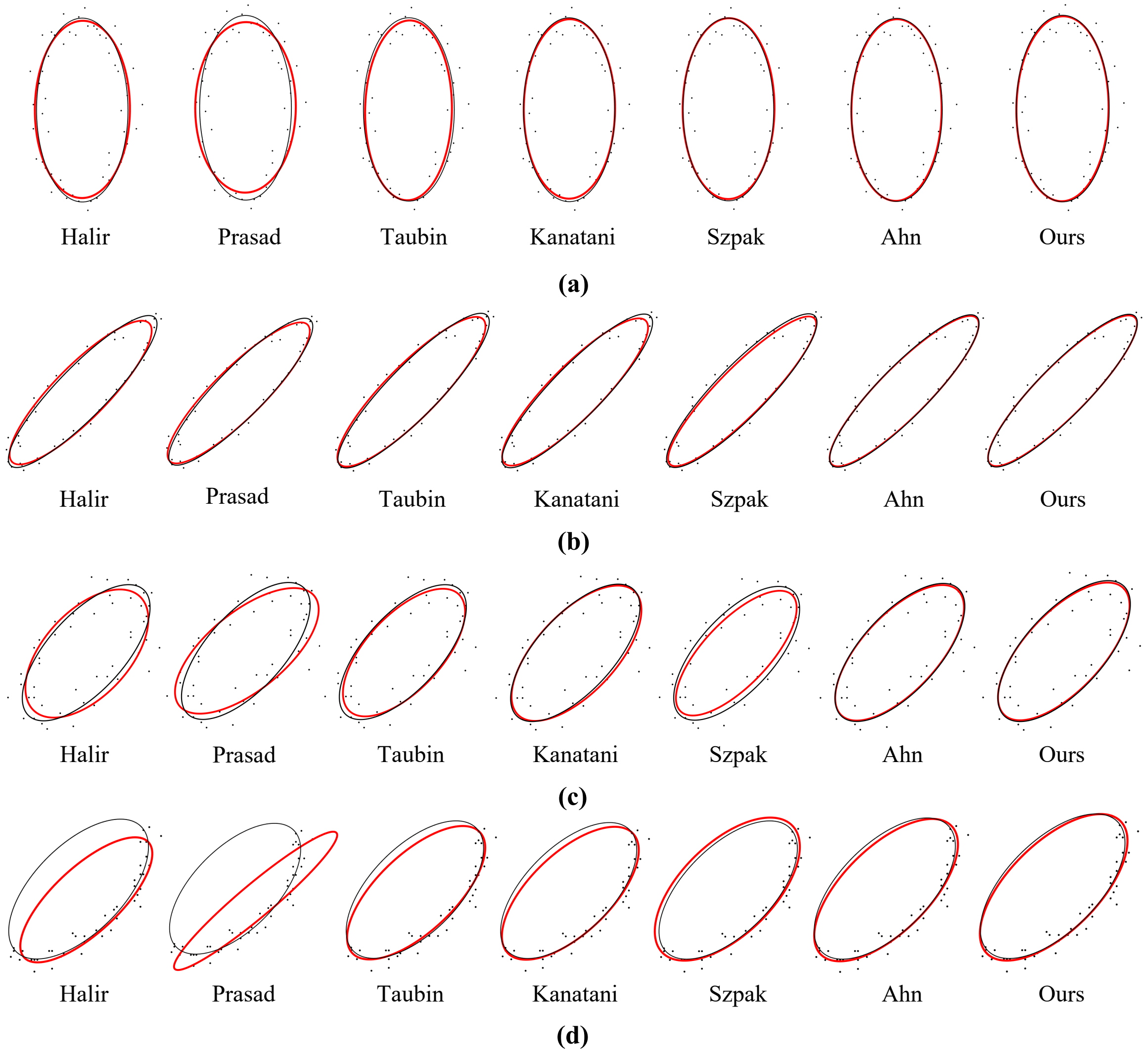

**Figure 11:** Example of the behavior of ellipse fitting methods (black ellipse is ground truth) of Halir, Prasad, Taubin, Kanatani, Szpak, Ahn and ours with the base configuration of Table 2 and the variable change: a) rotation angle of $\frac{\pi}{2}$ (vertical ellipse); b) aspect ratio of 4; c) noise of 4 pixels; and d) spanning arc angle of $\pi$ (half ellipse)

To summarize the results, the mean of the parameter errors obtained by each method is presented in Table 5 with the rows representing the results obtained by the methods for each experiment. As observed, for all experiments, ours and Ahn's gold standard geometric fit achieved nearly identical parameter estimation results, and outperformed the methods of Halir, Taubin, Kanatani and Szpak. It is worth noting that the spanning arc for the simulated ellipses in the first three experiments was $2\pi$ (full ellipse). When the arc angle was changed, however, our method even outperformed that of Ahn's (Table 5; row 4). The results presented in Tables 3 and 4 show that our method achieved consistently reliable ellipse parameter estimation results, comparable to established geometric fits (Ahn's method), and outperformed reliable algebraic (Halir) and Sampson-based methods (Taubin, Kanatani and Szpak). It was also observed that given the same initial ellipse estimate and convergence criteria, our method converged faster compared to the other two iterative methods (Ahn and Szpak), and consistently converged within 9-10 iterations on average. The results of the simulation experiments demonstrate that the confocal hyperbola distance is an excellent predictor of the true geometric distance and can also be used to produce reliable ellipse fitting results.

**Table 5:** Summary of simulated results: the mean of the parameter errors (%) of different methods

| Experiment | Halir | Taubin | Kanatani | Szpak | Ahn | Our Method |
|---|---|---|---|---|---|---|
| Rotation Angle | 0.450 | 0.430 | 0.419 | 0.412 | 0.399 | **0.397** |
| Aspect Ratio | 0.445 | 0.421 | 0.409 | 0.391 | 0.387 | **0.386** |
| Noise | 0.642 | 0.572 | 0.528 | 0.569 | 0.521 | **0.519** |
| Arc Angle | 7.170 | 1.917 | 1.881 | 1.739 | 1.736 | **1.589** |

Note: The bold numbers represent the best method for a given row

## 6. Real-world experiments

In this section, the results of fitting ellipses using the methods of Halir [16], Szpak [12], Ahn [27], and ours on 10,000 randomly oriented plane intersections with the cylinder point clouds of Figure 3 are reported. Table 6 shows the results of the mean of the parameter estimation errors for the center, radius, and axis angle of the cylinders of each case w.r.t. the ground truth. Four main observations are made from the results presented in Table 6. Firstly, the estimated center and radius of the cylinders appear to be directly impacted by the normalized RMSE ratio. In other words, the errors in the estimated center and radius (which represents the minor axis of the ellipses) appear to increase as the noise ratio increases. Second, the estimated center and radius using Halir's method is impacted the most (around 3mm and 2.5mm, respectively) with the increase of the normalized RMSE, compared to say our method (which changed around 1.5mm and 0.7mm, respectively). Based on our observations in Section 5.2.3, Halir's method started producing relatively biased results at noise levels above around 1.5 pixels noise (normalized RMSE of around 0.01); hence, it was expected that Halir's method produces more biased results compared to the remaining methods at larger noise levels in cases of Figures 3b and 3c. Halir's method also produced relatively more biased results compared to the remaining methods in each case, which was also expected since the method was shown to be produce considerably biased results for arc angles of less than around ~200° (see Figure 9). The third observation is that the axis angle error, which represents the ratio between the semi-minor and semi-major lengths, does not follow the same trend as the center and radius w.r.t. the relative noise increase. In fact, the ratio between the semi-minor and semi-major lengths appears to improve from case 3b to case 3c with very little change between cases 3a and 3c. Upon further investigation, it was observed that due to the lower point density of case 3b, the average number of points within 1mm of the intersecting planes was 30 in opposed to 45 and 48 points for cases 3a and 3c, respectively, which impacts the estimated ellipse parameters, including the ratio of the estimated semi-minor to semi-major lengths. The fourth observation is that even though Ahn's method provided more accurate radius estimations (corresponding to the semi-minor length) compared to Szpak's method, Szpak's method achieved lower axis angle errors, particularly in cases 3a and 3b. The difference is small (in the order of $10^{-1}$ degrees) and, hence, not a concern in practical settings. However, this shows that in these configurations, Szpak's method provides a slightly better estimation of the ratio between the semi-minor and semi-major lengths compared to Ahn's. The final observation from the data was that our method outperformed all other methods in the parameter estimation of the center, radius, and axis angle errors on average. On average for all measurements, our method achieved 1.5%, 5.7% and 85.3% better center estimation results, 1.7%, 6.0%, and 99.7% better radius estimation results, and 5.2%, 2.9% and 23.7% better cylinder axis estimation results, compared to the methods of Ahn, Szpak, and Halir, respectively. The average of the absolute difference in the estimation of the center, radius, and axis angle using our method were 0.2 mm, 0.2 mm and 0.1 degrees, respectively, compared to Szpak's, and 0.1 mm, 0.1 mm and 0.2 degrees, respectively, compared to Ahn's method, which are still notable in application pertaining to high precision metrology and reliability.

From the results of the simulation experiments presented in Section 5, the methods of Szpak, Ahn and ours behaved similarly in most configurations, other than for noise levels above 2.5 pixels (normalized RMSE of around 0.02 and above), where Szpak's method started to produce relatively biased results. The configurations in the real-world experiment were also in the normalized RMSE range, where Szpak's method was expected to behave similarly to ours and Ahn's. In larger noise levels, however, as shown in Figure 8, ours and Ahn's are expected to perform better. In Szpak [12], it was mentioned that the Sampson based methods are a good approximation for "moderate" noise levels. From our experiments, the 0.02 normalized RMSE ratio might be a good quantification of the definition of "moderate" in this context. As a point of reference, Figure 12 is provided as an example of the behavior of each ellipse fitting method with respect to the ground truth parameters for the projected points of Figure 3d-bottom. As illustrated from this example (Figure 12), our method appears to perform slightly better than Ahn, which is slightly better than Szpak, and considerably better than Halir.

**Table 6:** Results of the mean of the parameter estimation errors for the cases presented in Figure 3

| Mean of Cylinder Parameter Errors (Case of Figure 3a) | | | | |
|---|---|---|---|---|
| **Cylinder Parameter** | **Halir** | **Szpak** | **Ahn** | **Ours** |
| **Center (mm)** | 4.88 | 3.18 | 3.13 | **3.06** |
| **Radius (mm)** | 4.53 | 3.02 | 2.78 | **2.74** |
| **Axis Angle (°)** | 2.83 | 1.94 | 1.99 | **1.93** |
| Mean of Cylinder Parameter Errors (Case of Figure 3b) | | | | |
| **Cylinder Parameter** | **Halir** | **Szpak** | **Ahn** | **Ours** |
| **Center (mm)** | 7.13 | 3.58 | 3.27 | **3.20** |
| **Radius (mm)** | 7.03 | 3.15 | 3.08 | **3.01** |
| **Axis Angle (°)** | 3.95 | 4.15 | 4.34 | **3.87** |
| Mean of Cylinder Parameter Errors (Case of Figure 3c) | | | | |
| **Cylinder Parameter** | **Halir** | **Szpak** | **Ahn** | **Ours** |
| **Center (mm)** | 8.05 | 4.70 | **4.64** | **4.64** |
| **Radius (mm)** | 6.89 | 3.55 | 3.49 | **3.44** |
| **Axis Angle (°)** | 2.61 | 2.15 | 2.14 | **2.13** |

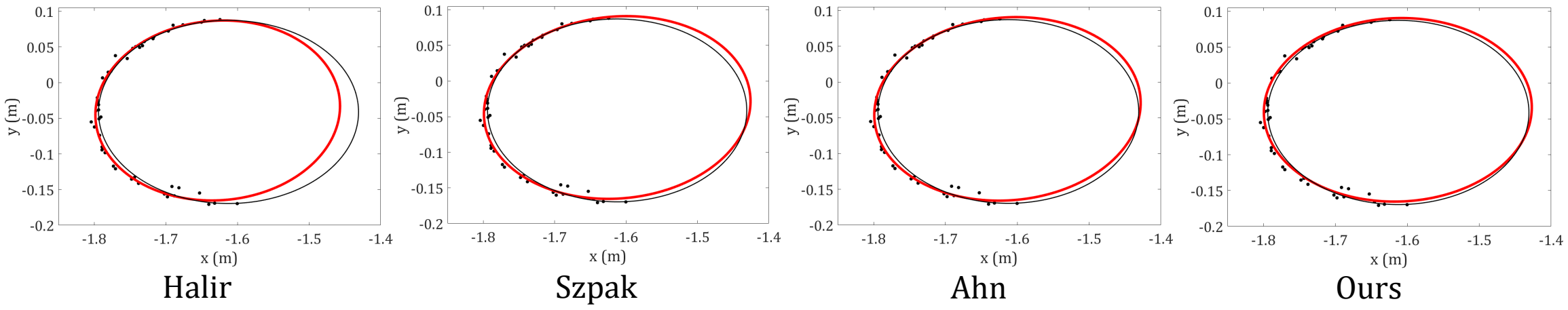

**Figure 12:** Behavior of the ellipse fitting methods of Halir, Szpak, Ahn and Ours for the projected cylindrical points of Figure 3d-bottom

## 7. Conclusions

This manuscript provided a new method for fitting ellipses to 2D points using the confocal hyperbola point to ellipse geometric distance approximation. Ellipse fitting using confocal hyperbola had not been formulated in the literature and, hence, its performance was never investigated. To this end, two comprehensive sets of experiments, including simulation-based and real-world were designed. The effectiveness of the proposed method compared to established state-of-the-art ellipse fitting methods in various configurations was examined. The simulation-based experiments consisted of two main categories, namely, point to ellipse distance approximations, and best fit ellipse evaluation. For point to ellipse geometric distance approximation methods, the confocal hyperbola method outperformed other established geometric distance approximations, namely, algebraic, Sampson, and Harker and O'Leary, and behaved almost identical (in the order of $10^{-2}$ pixels on average) to the ground truth geometric distance, proposed by Chernov [11] . To this end, a new ellipse fitting method using confocal hyperbola was proposed (Algorithm 1) and its effectiveness was evaluated compared to the ellipse fitting methods proposed by Halir [16], Taubin [19], Kanatani [20], Szpak [12], and Ahn [27]. The effects of change of variables, namely, rotation angle, aspect ratio, noise, and spanning arc on the estimated ellipse parameters using each method were thoroughly investigated. Overall, it was observed that the proposed confocal hyperbola-based ellipse fitting method achieved almost identical to (and in some cases outperformed) the gold standard geometric ellipse fit of Ahn, and outperformed the remaining methods in all categories. It was also shown that in all simulations, our method converged after around 9-10 iterations and achieved an average computation speed of roughly 140, 21 and 1.2 times faster than those achieved using the methods of Ahn, Szpak, and Kanatani.

The real-world experiment was designed to determine the ability of different ellipse fitting methods in predicting the parameters of cylindrical point clouds, acquired from mechanical pipes. To this end three cylindrical point clouds with different noise levels, and point densities were considered. The methods of Halir, Szpak, Ahn, and ours, which are all methods guaranteeing ellipses, were considered to estimate ellipse parameters of intersecting planes to these cylinders. The ellipse parameters were then used to recover the original cylinder parameter. The errors in estimating the original cylinder parameters were then recorded. It was observed that our method outperformed the remaining methods in estimating the cylinder's parameters. The method of Halir was clearly less accurate than the remaining methods, which was expected due to the biased behavior of Hair's ellipse parameter estimation in smaller arc angles (approximately less than a semi-ellipse) and at larger noise levels.

Overall, the comprehensive evaluation in both simulation-based and real-world experiments showed the agreement between the ellipse parameters estimated using the proposed confocal hyperbola ellipse fitting method and the true geometric parameters of the ellipses. The method achieved the best parameter estimation results in all categories of simulation and real-world experiments. We conclude that the confocal hyperbola-based ellipse fitting method shows great promise for fitting ellipses to 2D data in practical settings.

## Acknowledgement

The authors wish to acknowledge the support provided by the MJS Mechanical Ltd. and Michael Baytalan for their cooperation during data acquisition. This research project was partly funded by the Natural Sciences and Engineering Research Council (NSERC) of Canada (542980 - 19) and Alberta Innovates (G2020000051).

**Appendix I: Closed-form Jacobian matrix of the confocal hyperbola distance function**

This section describes the formulas to acquire the Jacobian of the distance function used to minimize equation 17 w.r.t. the geometric parameter vector $\rho = (x_c, y_c, a_e, b_e, \theta)$:

$$D_h = \sqrt{{D_{hX}}^2 + {D_{hY}}^2} \rightarrow \frac{\partial D_h}{\partial \rho} = \frac{D_{hX}.\frac{\partial D_{hX}}{\partial \rho} + D_{hY}.\frac{\partial D_{hY}}{\partial \rho}}{\sqrt{{D_{hX}}^2 + {D_{hY}}^2}}, \begin{cases} \frac{\partial D_{hX}}{\partial \rho} = \frac{\partial |X|}{\partial \rho} - \frac{\partial |X_I|}{\partial \rho} \\ \frac{\partial D_{hX}}{\partial \rho} = \frac{\partial |Y|}{\partial \rho} - \frac{\partial |Y_I|}{\partial \rho} \end{cases} \tag{A1}$$

$$|X_I| = \frac{a_e.|X|}{\sqrt{\frac{T+\sqrt{\Delta}}{2}}}, \text{where:} \begin{cases} T = X^2 + Y^2 + f^2 \\ \Delta = T^2 - 4X^2 f^2 \end{cases} \rightarrow \frac{\partial |X_I|}{\partial \rho} = \frac{\frac{\partial |X|}{\partial \rho} a_e}{\sqrt{\frac{T+\sqrt{\Delta}}{2}}} + \frac{\frac{\partial a_e}{\partial \rho}.|X|}{\sqrt{\frac{T+\sqrt{\Delta}}{2}}} - \frac{a_e.|X|.(2\sqrt{\Delta}\frac{\partial T}{\partial \rho} + \frac{\partial \Delta}{\partial \rho})}{8.\sqrt{\Delta\left(\frac{T+\sqrt{\Delta}}{2}\right)^3}} \tag{A2}$$

$$|Y_I| = b_e\sqrt{1-L^2}, \text{where: } L = \frac{X_I}{a_e} = \frac{|X|}{\sqrt{\frac{T+\sqrt{\Delta}}{2}}} \rightarrow \frac{\partial |Y_I|}{\partial \rho} = \frac{\partial b_e}{\partial \rho}\sqrt{1-L^2} - \frac{b_e.L.\left(\frac{\frac{\partial |X|}{\partial \rho}}{\sqrt{\frac{T+\sqrt{\Delta}}{2}}} - \frac{|X|.\left(2\sqrt{\Delta}\frac{\partial T}{\partial \rho} + \frac{\partial \Delta}{\partial \rho}\right)}{8\sqrt{\Delta\left(\frac{T+\sqrt{\Delta}}{2}\right)^3}}\right)}{\sqrt{1-L^2}} \tag{A3}$$

$$\begin{cases} \frac{\partial a_e}{\partial \rho} = (0,0,1,0,0) \\ \frac{\partial b_e}{\partial \rho} = (0,0,0,1,0) \\ \frac{\partial |X|}{\partial \rho} = \text{sign}(X).\frac{\partial X}{\partial \rho} = \text{sign}(X).(-C,-S,0,0,Y) \\ \frac{\partial |Y|}{\partial \rho} = \text{sign}(Y).\frac{\partial Y}{\partial \rho} = \text{sign}(Y).(S,-C,0,0,-X) \\ \frac{\partial T}{\partial \rho} = 2X\frac{\partial X}{\partial \rho} + 2Y\frac{\partial Y}{\partial \rho} + 2a_e\frac{\partial a_e}{\partial \rho} - 2b_e\frac{\partial b_e}{\partial \rho} \\ \frac{\partial \Delta}{\partial \rho} = 2T\frac{\partial T}{\partial \rho} - 8X^2(a_e\frac{\partial a_e}{\partial \rho} - b_e\frac{\partial b_e}{\partial \rho}) - 8f^2 X\frac{\partial X}{\partial \rho} \end{cases} \tag{A4}$$

where $C = \cos\theta$, $S = \sin\theta$, $X = (x - x_c).C + (y - y_c).S$, and $Y = -(x - x_c).S + (y - y_c).C$. From the above formulations, three possible numerical singularities in the Jacobian, $JD_h$, may appear when either: (i) $L = 1$ (equation A3); (ii) $\Delta = 0$ (equation A2); or (iii) when $D_h = 0$ (equation A1). The case of $L = 1$ belongs to the situation when the intersecting point is on the semi-major length of the ellipse i.e. when $|X| \geq f$ and $|Y| = 0$. The second case $\Delta = 0$ occurs when $|X| = f$ and $|Y| = 0$, which is a condition covered more broadly in the previous case. These cases occur when $|X_I| = a_e$. For case 1 (and consequentially 2), equation 16 provided a simplified version of the confocal hyperbola distance function (see Section 3), which is recommended to be used here. The Jacobian matrix for the conditions of equation 16 (covering cases 1 and 2) is provided in equation A6 below. The final possible numerical singularity might occur when the points are exactly on the ellipse (i.e. $D_{hX} = D_{hY} = 0$). For this situation, equation A1 can be re-written and the Jacobian can be approximated by taking its limit as $(D_{hX}, D_{hY}) \rightarrow (0,0)$:

$$
\begin{cases}
\dfrac{\partial D_h}{\partial \rho} = \dfrac{\frac{\partial D_{hX}}{\partial \rho}}{\sqrt{1+\left(\frac{D_{hY}}{D_{hX}}\right)^2}} + \dfrac{\frac{\partial D_{hY}}{\partial \rho}}{\sqrt{\left(\frac{D_{hX}}{D_{hY}}\right)^2+1}} \\
\lim\limits_{(D_{hX},D_{hY})\to(0,0)} \dfrac{\partial D_h}{\partial \rho} \equiv \lim\limits_{(D_{hX},D_{hX})\to(0,0)} \dfrac{\partial D_h}{\partial \rho} = \dfrac{\frac{\partial D_{hX}}{\partial \rho}}{\sqrt{1+\left(\frac{D_{hX}}{D_{hX}}\right)^2}} + \dfrac{\frac{\partial D_{hY}}{\partial \rho}}{\sqrt{\left(\frac{D_{hX}}{D_{hX}}\right)^2+1}} = \dfrac{\sqrt{2}}{2}\left(\dfrac{\partial D_{hX}}{\partial \rho} + \dfrac{\partial D_{hY}}{\partial \rho}\right)
\end{cases}
\tag{A5}
$$

Equation A5 is now numerically stable for points on the ellipse and can also be used for any other point whose distances in $X$ and $Y$ are equal (i.e. $D_{hX} = D_{hY}$). These formulations can be summarized as follows:

$$
\begin{bmatrix} \frac{\partial D_{hX}}{\partial \rho} \\ \frac{\partial D_{hY}}{\partial \rho} \end{bmatrix} =
\begin{cases}
\begin{bmatrix} 0 & 0 & 0 & 0 & 0 \\ \text{sign}(Y).S & -\text{sign}(Y).C & 0 & -1 & -\text{sign}(Y).X \end{bmatrix} & , \quad |Y_I| = b_e \\
\begin{bmatrix} -\text{sign}(X).C & -\text{sign}(X).S & -1 & 0 & \text{sign}(X).Y \\ 0 & 0 & 0 & 0 & 0 \end{bmatrix} & , \quad |X_I| = a_e \\
\text{Use equations } A1 - A4 & , \quad \text{Anywhere else}
\end{cases}
\tag{A6}
$$

$$
\frac{\partial D_h}{\partial \rho} =
\begin{cases}
\dfrac{D_{hX}.\frac{\partial D_{hX}}{\partial \rho} + D_{hY}.\frac{\partial D_{hY}}{\partial \rho}}{\sqrt{{D_{hX}}^2 + {D_{hY}}^2}} & , \quad D_{hX} = 0 \text{ or } D_{hY} = 0 \\
\dfrac{\sqrt{2}}{2}\left(\dfrac{\partial D_{hX}}{\partial \rho} + \dfrac{\partial D_{hY}}{\partial \rho}\right) & , \quad D_{hX} = D_{hY} \\
\dfrac{\frac{\partial D_{hX}}{\partial \rho}}{\sqrt{1+\left(\frac{D_{hY}}{D_{hX}}\right)^2}} + \dfrac{\frac{\partial D_{hY}}{\partial \rho}}{\sqrt{\left(\frac{D_{hX}}{D_{hY}}\right)^2+1}} & , \quad \text{Anywhere else}
\end{cases}
\tag{A7}
$$